\documentclass[preprint,10pt]{elsarticle}

\usepackage{amssymb}

\usepackage{amsmath,amsthm}
\numberwithin{equation}{section}
\usepackage{dsfont}
 \usepackage{amsopn}
\DeclareMathOperator{\Supp}{Supp}
\DeclareMathOperator{\Div}{div}

\theoremstyle{definition} 

\theoremstyle{remark}
\newtheorem{rem}{\textsf{Remark}}

\theoremstyle{plain}
\newtheorem{theo}{\textsc{Theorem}}
\newtheorem{Lemme}{\textsf{Lemma}}[section]
\newtheorem{prop}{\textsc{Proposition}}[section]

\newcommand{\cqfd}
{\mbox{}\nolinebreak\hfill\rule{2mm}{2mm}\medbreak\par}

\usepackage{geometry}
\journal{Journal de Math\'{e}matiques Pures et Appliqu\'{e}es}

\begin{document}

\begin{frontmatter}



\title{Indirect controllability of some linear parabolic systems of $m$ equations with 
$m-1$ controls involving coupling terms of zero or first order.}


 \author[label1]{Michel Duprez}
 \author[label2]{Pierre Lissy}
 \address[label1]{Laboratoire de Math\'ematiques de Besan\c{c}on, UMR CNRS 6623, Universit\'e de Franche-Comt\'e,
16 route de Gray, 25030 Besan\c{c}on Cedex, France, \texttt{michel.duprez@univ-fcomte.fr}}
 \address[label2]{CEREMADE, UMR CNRS 7534, Universit\'e Paris-Dauphine, Place du Mar\'echal de Lattre de Tassigny, 75775 Paris Cedex
16, France, \texttt{lissy@ceremade.dauphine.fr.}}

\begin{abstract}
 This paper is devoted to the study of the null and approximate controllability 
 for some classes of linear 
 coupled parabolic systems with less controls than equations. More precisely, 
 for a given bounded domain $\Omega$ 
 in $\mathbb{R}^N$ ($N\in \mathbb{N}^*$), we consider a system of $m$ 
 linear parabolic equations ($m\geqslant 2$)  with coupling terms of first and zero order, and $m-1$ controls 
 localized in some arbitrary nonempty open subset $\omega$ of $\Omega$. In the case of constant coupling 
 coefficients, we provide a necessary and sufficient condition to obtain the null or approximate
 controllability in arbitrary small time. In the case $m=2$ and $N=1$, 
 we also give a generic sufficient condition   to obtain the null or approximate controllability in arbitrary small time for general coefficients 
 depending on the space and times variables, provided that the supports of the coupling terms intersect the control domain $\omega$. 
The results are obtained thanks to the fictitious control method together with an algebraic 
method  and some appropriate Carleman estimates.
\end{abstract}


\end{frontmatter}

\noindent{\bf R\'esum\'e}\\
Cet article est consacr\'e \`a l'\'etude de la contr\^olabilit\'e \`a z\'ero et approch\'ee 
d'une classe de syst\`emes paraboliques lin\'eaires coupl\'es avec moins de contr\^oles 
que d'\'equations. 
Plus pr\'ecis\'ement, pour un domaine donn\'e born\'e $\Omega$ de $\mathbb{R}^N$ 
($N\in\mathbb{N}^*$), nous consid\'erons un syst\` eme de $m$ \'equations ($m\geqslant2$) 
avec des termes de couplages d'ordre un et z\'ero, et $m-1$ contr\^oles localis\'es dans un 
sous-ensemble ouvert non-vide arbitraire $\omega$ de $\Omega$. 
Dans le cas de coefficients de couplage constants, 
nous fournissons une condition n\'ecessaire et suffisante pour obtenir la contr\^olabilit\'e 
\`a z\'ero ou approch\'ee en temps arbitrairement petit.
Dans le cas $m=2$ et $N=1$, nous donnons \'egalement une condition g\'en\'erique suffisante 
pour obtenir la contr\^olabilit\'e \`a z\'ero ou approch\'ee en temps arbitrairement petit 
pour des coefficients g\'en\'eraux d\'ependants des variables de temps et d'espace, 
lorsque le support des termes de couplage intersecte le domaine de contr\^ole $\omega$. 
Les r\'esultats sont obtenus par combinaison de la m\'ethode par contr\^ole fictif 
avec une m\'ethode alg\'ebrique et des estimations de Carleman appropri\'ees.





\section{Introduction}
\subsection{Presentation of the problem and main results}

 Let $T>0$, let $\Omega$ be a bounded domain in $\mathbb{R}^N$ ($N\in\mathbb{N}^*$) supposed to be regular enough (for example of class $\mathcal{C}^{\infty}$),
and let $\omega$ be an arbitrary nonempty open subset of $\Omega$. Let $Q_T:=(0,T)\times\Omega$, 
$q_T:=(0,T)\times\omega$, $\Sigma_T:=(0,T)\times\partial\Omega$ and  $m\geqslant 2$. 
We consider the following system of $m$ parabolic linear equations, where the $m-1$ first equations are controlled:
\begin{equation}\label{system primmal}
 \left\{\begin{array}{ll}
\partial_ty_1=\Div(d_1\nabla y_1)+\sum_{i=1}^mg_{1i}\cdot\nabla y_i+\sum_{i=1}^{m}a_{1i}y_i+\mathds{1}_{\omega}u_1&\mathrm{in}~ Q_T,\\
\partial_ty_2=\Div(d_2\nabla y_2)+\sum_{i=1}^mg_{2i}\cdot\nabla y_i+\sum_{i=1}^{m}a_{2i}y_i+\mathds{1}_{\omega}u_2&\mathrm{in}~ Q_T,\\
\vdots\\
\partial_ty_{m-1}=\Div(d_{m-1}\nabla y_{m-1})+\sum_{i=1}^mg_{(m-1)i}\cdot\nabla y_i+\sum_{i=1}^{m}a_{(m-1)i}y_i+\mathds{1}_{\omega}u_{m-1}&\mathrm{in}~ Q_T,\\
\partial_ty_m=\Div(d_m\nabla y_m)+\sum_{i=1}^mg_{mi}\cdot\nabla y_i+\sum_{i=1}^{m}a_{mi}y_i&\mathrm{in}~ Q_T,\\
y_1=\ldots=y_m=0&\mathrm{on}~\Sigma_T,\\
y_1(0,\cdot)=y_1^0,\ldots,~y_m(0,\cdot)=y_m^0&\mathrm{in}~\Omega,
        \end{array}
\right.
\end{equation}
where $y^0:=(y_1^0,\ldots,y_m^0)\in L^2(\Omega)^m$ is the initial condition 
and $u:=(u_1,\ldots,u_{m-1})\in L^2(Q_T)^{m-1}$ is the control. 
The  zero and first  order coupling terms
  $(a_{ij})_{1\leqslant i,j\leqslant m}$ and $(g_{ij})_{1\leqslant i,j\leqslant m}$ 
  are assumed to be respectively in  $ L^\infty(Q_T)$ and 
 in $L^\infty(0,T;W^{1}_{\infty}(\Omega)^N)$.
Given some $l\in\{1,...,m\}$, the second order elliptic self-adjoint operator $\Div(d_l\nabla)$ is given by
 \begin{equation*}
  \Div(d_l\nabla)=\sum\limits_{i,j=1}^N\partial_i(d^{ij}_l\partial_j),
 \end{equation*}
with 
\begin{equation*}\left\{
 \begin{array}{l}
  d^{ij}_l\in W^{1}_{\infty}(Q_T),\\
   d^{ij}_l= d^{ji}_l\mathrm{~in~}Q_T,
 \end{array}\right.
\end{equation*}
where the coefficients $d^{ij}_l$ satisfy the uniform ellipticity condition 
\begin{equation*}
\sum\limits_{i,j=1}^N d^{ij}_l\xi_i\xi_j\geqslant d_0|\xi|^2\mathrm{~in~}Q_T,~\forall \xi\in\mathbb{R}^N,
\end{equation*}
for a constant $d_0>0$.

In order to simplify the notation, from now on, we will denote by
\begin{equation*}
\begin{aligned}
 D:=\mathrm{diag}(d_1,\ldots, d_m),~~&
G:=(g_{ij})_{1\leqslant i,j\leqslant m}\in \mathcal M_m(\mathbb R^N),\\
  A:=(a_{ij})_{1\leqslant i,j\leqslant m}\in \mathcal M_m(\mathbb R), ~&
  B:=\left(\begin{array}{c}\mathrm{diag}(1,\ldots,1)\\0\end{array}\right)\in \mathcal M_{m,m-1}(\mathbb R),
  \end{aligned}
\end{equation*}
so that we can write System \eqref{system primmal} as 
\begin{equation}\label{system primmal2}
 \left\{\begin{array}{ll}
\partial_ty=\Div(D\nabla y)+G\cdot\nabla y+Ay+\mathds{1}_{\omega}Bu&\mbox{ in } Q_T,\\
y=0&\mbox{ on }\Sigma_T,\\
y(0,\cdot)=y^0&\mbox{ in }~\Omega.
        \end{array}
\right.
\end{equation}

It is well-known (see for instance \cite[Th. 3 \& 4, p. 356-358]{MR2597943}) 
that for  any initial data  $y^0\in L^2(\Omega)^m$ and $u\in L^2(Q_T)^{m-1}$, 
 System \eqref{system primmal2} 
 admits a unique solution $y$  
in 
$W(0,T)^m$, where  
\begin{equation}\label{defW}
W(0,T):=L^2(0,T;H^1_0(\Omega ))
\cap H^1(0,T;H^{-1}(\Omega ))\hookrightarrow\mathcal{C}^0([0,T];L^2(\Omega)).
\end{equation}
Moreover, one can prove (see for instance \cite[Th. 5, p. 360]{MR2597943}) that if $y^0\in H^1_0(\Omega)^m$ and  $u\in L^2(Q_T)^{m-1}$, then
the solution $y$  is in 
$W^{2,1}_2(Q_T)^m$, where  
\begin{equation}\label{defW21}
W^{2,1}_2(Q_T):= L^2(0,T;H^2(\Omega)\cap H^1_0(\Omega ))
 \cap H^1(0,T;L^2(\Omega ))\hookrightarrow\mathcal{C}^0([0,T];H^1_0(\Omega)).
\end{equation}

The main goal of this article is to analyse the null controllability and approximate controllability
of System \eqref{system primmal}. Let us recall the definition of these notions. It will be said that
\begin{itemize}
 \item[$\bullet$] System \eqref{system primmal} is \textit{null controllable} at time $T$  
if for every initial condition  $y^0\in  L^2(\Omega)^m$, there exists a control $u\in L^2(Q_T)^{m-1}$ 
such that the solution $y$ in $W(0,T)^m$ to System \eqref{system primmal} satisfies
\begin{equation*}
 y(T)\equiv0\mathrm{~~in~~}\Omega.
\end{equation*}
 \item[$\bullet$] System \eqref{system primmal} is \textit{approximately controllable} 
 at time $T$
if for every $\varepsilon>0$, every initial condition  $y^0\in  L^2(\Omega)^m$ and 
every $y_T\in  L^2(\Omega)^m$, there exists a control $u\in L^2(Q_T)^{m-1}$ 
such that the solution $y$ in  $W(0,T)^m$ to System \eqref{system primmal} satisfies
\begin{equation*}
 \|y(T)-y_T\|^2_{L^2(\Omega)^m}\leqslant \varepsilon.
\end{equation*}
\end{itemize}

 Let us remark that if  System \eqref{system primmal} is null controllable on the time interval $(0,T)$, 
 then it is also approximately controllable on the time interval $(0,T)$
 (this is an easy consequence of usual results 
 of backward uniqueness concerning parabolic equations as given for example in \cite{MR0338517}).


Our first  result  gives a necessary and sufficient condition 
for null (or approximate)  controllability of System \eqref{system primmal} in the case of constant coefficients. 

\begin{theo}\label{th 1}
Let us assume that $D$, $G$ and $A$ are \textbf{\emph{constant in space and time}}. 
Then System \eqref{system primmal} is null (resp. approximately) controllable at time $T>0$  
 if and only if there exists $i_0\in\{1,...,m-1\}$ such that 
 \begin{equation}\label{cond th 1}\begin{array}{c}
g_{mi_0}\neq0
  ~~\mathrm{or}~~
a_{mi_0}\neq0.
 \end{array}\end{equation}
\end{theo}
This condition is the natural one that can be expected, since it means that the last equation is coupled with at least one of the others, and is clearly a necessary condition (otherwise one cannot act on the last component $y_n$ which evolves freely).

\vspace{0,2cm}

Our second result concerns the case of general coefficients depending on  space and time variables, 
in the particular case of two equations (i.e. $m=2$), 
and gives a controllability result under some technical conditions
 on the coefficients (see \eqref{cond1th2} and \eqref{cond2th2}) coming from the algebraic solvability 
 (see Section $3.1$). To understand why this kind of condition appears here, 
 we refer to the simple example given in \cite[Ex. 1, Sec. 1.3]{thL}. Let us emphasize that the second point is only valid for $N=1$ and under Condition \eqref{cond2th2}, which is clearly technical since it does not even cover the case of constant coefficients. However, Condition \eqref{cond2th2} is \emph{generic} as soon as we restrict to the coupling coefficients verifying $g_{21}\not =0$ on $(0,T)\times\omega$, in the following sense: it only requires some regularity on the coefficients 
 and a given determinant, involving some coefficients and their derivatives, to be  non-zero. 
 Since this condition may seem a little bit intricate, we will give in Remark~\ref{rem1} some particular examples that clarify the scope of Theorem $2$.

\begin{theo}\label{th 2}
Consider the following system:
\begin{equation}\label{syst21}
 \left\{\begin{array}{ll}
\partial_ty_1=\Div(d_1\nabla y_1)+g_{11}\cdot\nabla y_1+g_{12}\cdot\nabla y_2+a_{11}y_1+a_{12}y_2
+\mathds{1}_{\omega}u&\mathrm{in}~Q_T,\\
\partial_ty_2=\Div(d_2\nabla y_2)+g_{21}\cdot\nabla y_1
+g_{22}\cdot\nabla y_2+a_{21}y_1+a_{22}y_2&\mathrm{in}~ Q_T,\\
y_1=y_2=0&\mathrm{on}~\Sigma_T,\\
y_1(0,\cdot)=y_1^0,~y_2(0,\cdot)=y_2^0,&\mathrm{in}~\Omega,
        \end{array}
\right.
\end{equation}
where $y^0=(y^0_1,y^0_2)\in L^2(\Omega)^2$ is the initial condition.

Then System \eqref{system primmal} is null (resp. approximately)  controllable  at time $T$  
 if there exists an open subset $(a,b)\times{\mathcal O}\subseteq q_T$ where 
 one of the following conditions is verified:
 \begin{enumerate}
  \item[(i)] Coefficients of System \eqref{syst21} satisfy  
  $d_i\in \mathcal C^1((a,b),\mathcal C^2(\mathcal O)^{N^2})$, 
    $g_{ij}\in \mathcal C^1((a,b),\mathcal C^2(\mathcal O)^N)$, 
    $a_{ij}\in \mathcal C^1((a,b),\mathcal C^2(\mathcal O))$ for $i,j=1,2$ 
  and 
  \begin{equation}\label{cond1th2}
\left .
\begin{array}{ll}
    g_{21}=0 \mathrm{ ~and~ }
  a_{21}\not =0&\mathrm{~ in~ }(a,b)\times{\mathcal O}.\\
  \end{array}\right.
  \end{equation}
  \item[(ii)]$N=1$ and 
  coefficients of System \eqref{syst21} satisfy $d_i,~g_{ij},~a_{ij}\in \mathcal C^3((a,b),\mathcal C^7(\mathcal O))$ for $i=1,2$ and
  \begin{equation}\label{cond2th2}
\begin{array}{l}
  |\mathrm{det}(H(t,x))|>C  \mbox{ for every }(t,x)\in(a,b)\times{\mathcal O},
  \end{array}\end{equation}
 where
  \begingroup\footnotesize
  \begin{equation}\label{defH}
  H:=\left(\begin{array}{ccccccc} 
-a_{21}+\partial_x g_{21}&g_{21}&0&0&0&0\\
-\partial_x a_{21}+\partial_{xx} g_{21}&-a_{21}+2\partial_x g_{21}&0&g_{21}&0&0\\
-\partial_t a_{21}+\partial_{tx} g_{21}&\partial_{t} g_{21}&-a_{21}+\partial_x g_{21}&0&g_{21}&0\\
-\partial_{xx} a_{21}+\partial_{xxx} g_{21}&-2\partial_x a_{21}+3\partial_{xx} g_{21}&0&
- a_{21}+3\partial_{x} g_{21}&0&g_{21}\\
-a_{22}+\partial_xg_{22}&g_{22}-\partial_xd_2&-1&-d_2&0&0\\
-\partial_x a_{22}+\partial_{xx}g_{22}&-a_{22}+2\partial_xg_{22}-\partial_{xx}d_2&0&g_{22}
-2\partial_xd_2&-1&-d_2\end{array}\right ).
  \end{equation}\endgroup

 \end{enumerate}

\end{theo}
\begin{rem}\label{rem1}
\begin{enumerate}
\item[(a)] The first part of Theorem \ref{th 2} has already been proved in \cite{gonzalez2010controllability} and \cite{gonzalezperez2006} (see the point $4.$ of Section $9$ in  \cite{gonzalezperez2006}) with less regularity on the coefficients, and is not a new result.
\item[(b)] We will see during the proof of Theorem \ref{th 2} that, in Item (ii) of Theorem \ref{th 2}, taking into account the derivatives of the appearing in \eqref{defH}, \eqref{M*21} and \eqref{def S alg} 
and the regularity needed for the control in Proposition \ref{prop cont reg2}, 
we only need the following regularity for the coefficients:
    \begin{equation}\label{reg1th2}
   \begin{array}{l}
  d_{i}\in \mathcal C^1((a,b),\mathcal C^2(\mathcal O)),~
  g_{ij}\in \mathcal C^0((a,b),\mathcal C^2(\mathcal O)),~
  a_{ij}\in \mathcal C^0((a,b),\mathcal C^1(\mathcal O))
  \end{array}\end{equation}
for all $i,j\in\{1,2\}$ in the first case (i) and 
      \begin{equation}\label{reg2th2}
       \left\{\begin{array}{l}
       d_{1}\in \mathcal C^1((a,b),\mathcal C^4(\mathcal O))\cap \mathcal C^2((a,b),\mathcal C^1(\mathcal O)),\\
g_{11},\mbox{ }g_{12}\in \mathcal C^0((a,b),\mathcal C^4(\mathcal O))\cap \mathcal C^1((a,b),\mathcal C^1(\mathcal O)),\\
  a_{11},\mbox{ }a_{12}\in \mathcal C^0((a,b),\mathcal C^3(\mathcal O))\cap \mathcal C^1((a,b),\mathcal C^0(\mathcal O)),\\
   g_{21}\in \mathcal C^0((a,b),\mathcal C^7(\mathcal O))\cap \mathcal C^3((a,b),\mathcal C^1(\mathcal O)),\\
  d_{2},\mbox{ }g_{22}\in\mathcal C^0((a,b),\mathcal C^6(\mathcal O))\cap \mathcal C^2((a,b),\mathcal C^2(\mathcal O)),\\
    a_{22}\in\mathcal C^0((a,b),\mathcal C^5(\mathcal O))\cap \mathcal C^2((a,b),\mathcal C^1(\mathcal O)),\\
        a_{21}\in\mathcal C^0((a,b),\mathcal C^6(\mathcal O))\cap \mathcal C^3((a,b),\mathcal C^0(\mathcal O)),\\
  \end{array}\right.\end{equation}
in the second case $(ii)$.

\item[(c)]  One can easily compute explicitly the determinant of matrix $H$ appearing in \eqref{cond2th2}:
   \begingroup\scriptsize
  \begin{equation*}
 \begin{array}{l}\mathrm{det}(H)=
  2 \frac{\partial a_{21}}{\partial x} \frac{\partial d_{2}}{\partial x} g_{21}^2
  -4 \frac{\partial a_{21}}{\partial x} {d_{2}}
      \frac{\partial g_{21}}{\partial x} {g_{21}}+  
    \frac{\partial^2 a_{21}}{\partial x^2} {d_{2}} g_{21}^2+2 {a_{21}}
     \frac{\partial a_{21}}{\partial x} {d_{2}} {g_{21}}-  \frac{\partial a_{21}}{\partial x} g_{21}^2
    {g_{22}}+  \frac{\partial a_{21}}{\partial t} g_{21}^2\\-4 {a_{21}}  
    \frac{\partial d_{2}}{\partial x}
      \frac{\partial g_{21}}{\partial x} {g_{21}}+{a_{21}}   
      \frac{\partial^2 d_{2}}{\partial x^2} g_{21}^2+a_{21}^2
     \frac{\partial d_{2}}{\partial x} {g_{21}}-3 {a_{21}} {d_{2}}  
     \frac{\partial^2 g_{21}}{\partial x^2} {g_{21}}+6
    {a_{21}} {d_{2}}   (\frac{\partial g_{21}}{\partial x})^2-2 a_{21}^2 {d_{2}}
     \frac{\partial g_{21}}{\partial x}+{a_{21}} 
     \frac{\partial g_{21}}{\partial x} {g_{21}} {g_{22}}\\-{a_{21}}
     \frac{\partial g_{21}}{\partial t} {g_{21}}-{a_{21}} g_{21}^2
   \frac{\partial g_{22}}{\partial x}+  \frac{\partial a_{22}}{\partial x} g_{21}^3-  
   \frac{\partial^2 d_{2}}{\partial x^2}   \frac{\partial g_{21}}{\partial x}
    g_{21}^2-2   \frac{\partial d_{2}}{\partial x}  
    \frac{\partial^2 g_{21}}{\partial x^2} g_{21}^2+3  \frac{\partial d_{2}}{\partial x}
      (\frac{\partial g_{21}}{\partial x})^2 {g_{21}}-{d_{2}}  
      \frac{\partial^3 g_{21}}{\partial x^3} g_{21}^2\\+5
    {d_{2}}   \frac{\partial g_{21}}{\partial x}  
    \frac{\partial^2 g_{21}}{\partial x^2} {g_{21}}-4 {d_{2}}
     (\frac{\partial g_{21}}{\partial x})^3+  \frac{\partial g_{21}}{\partial x} g_{21}^2
   \frac{\partial g_{22}}{\partial x}+  \frac{\partial g_{21}}{\partial x^2} g_{21}^2 {g_{22}} 
   -  \frac{\partial g_{21}}{\partial x}^2
    {g_{21}} {g_{22}}-  \frac{\partial g_{21}}{\partial x\partial t} g_{21}^2
    +  \frac{\partial g_{21}}{\partial x}
     \frac{\partial g_{21}}{\partial t} {g_{21}}\\-g_{21}^3   \frac{\partial g_{22}}{\partial x^2}.
     \end{array}
  \end{equation*}
  \endgroup
  \item[(d)]
We remark that Condition \eqref{cond2th2} implies in particular that
\begin{equation}\label{g21 non nul}
 g_{21}\neq0\mathrm{~in~}(a,b)\times{\mathcal O}.
\end{equation}
Our conjecture is that, as in the case of constant coefficients, either the first line 
of \eqref{cond1th2} or \eqref{g21 non nul} is sufficient as soon as we restrict 
to the class of coupling terms that intersect the control region, since it is the minimal conditions one can expect (as in the case of constant coefficients, this only means that the last equation is coupled with one of the others).
\item[(e)]
Even though Condition \eqref{cond2th2} seems complicated, it can be simplified in some cases. 
Indeed, for example, System \eqref{syst21} is null controllable  at time $T$ 
if there exists an open subset 
$(a,b)\times{\mathcal O}\subseteq q_T$ such that 
  \begin{equation*}
\left\{\begin{array}{ll}
   g_{21}\equiv \kappa\in \mathbb{R}^*&\mathrm{~ in~ }(a,b)\times{\mathcal O},\\
 a_{21}\equiv0&\mathrm{~ in~ }(a,b)\times{\mathcal O},\\
  \partial_x a_{22}\neq \partial_{xx} g_{22}&\mathrm{~ in~ }(a,b)\times{\mathcal O}.
  \end{array}\right.
  \end{equation*}
In the case $\partial_x a_{22}= \partial_{xx} g_{22}$, we do not know if the controllability holds and we are not able to prove it using the same techniques as in this paper.

Another simple situation is the case where the coefficients depend only on the time variable. 
In this case, it is easy to check that Condition \eqref{cond2th2} becomes simply: 
there exists an open interval $(a,b)$ of $(0,T)$ such that
  \begin{equation*}
\left .\begin{array}{ll}
  g_{21}(t)\partial_t a_{21}(t)\neq a_{21}(t) \partial_t g_{21}(t)&\mathrm{~ in~ }(a,b).
  \end{array}\right.
  \end{equation*}

            \item[(f)] In fact, it is likely that one could obtain a far more general result 
than the one obtained in Theorem \ref{th 2}. By using the same reasoning, 
one would be able to obtain a result of controllability for arbitrary $m$ and $N$, 
however the generic Condition \eqref{cond2th2} would be far more complicated 
and in general impossible to write down explicitly. 
That is the reason why we chose to treat only the case $m=2$ and $N=1$.
\end{enumerate}
\end{rem}

\indent This paper is organized as follows. In Section $1.2$,  
we  recall some previous results and explain precisely the scope of the present contribution. In Section $1.3$, we give some natural perspectives of this work.
In Section $1.4$ we present the main method used here, that is to say the fictitious 
control method together with some algebraic method. 
Section 2 is devoted to the proof of Theorem \ref{th 1}.
We finish with the proof of Theorem \ref{th 2} in Section $3$.

\subsection{State of the art}

 The study of what is  called the \emph{indirect} controllability for linear or non-linear 
parabolic coupled systems has been an intensive subject of interest these last years. 
The main issue is to try to control many equations with less controls than equations 
(and ideally only one control if possible), with the hope that one can act \emph{indirectly} on 
the equations that are not directly controlled thanks to the coupling terms.
For a recent survey concerning  this kind of control problems, we refer to \cite{ammar2011recent}. 
Here, we will mainly present the results related to this work, that is to say the case 
of the null or approximate controllability  of linear parabolic systems with distributed controls.

  First of all, in the case of zero order coupling terms, 
   a necessary and sufficient algebraic condition is proved in \cite{ammar2009generalization}
  for the controllability of parabolic systems, for constant coefficients 
  and diffusion coefficients $d_i$ that are equal. 
  This condition is similar to the usual algebraic 
  \emph{Kalman rank condition} for finite-dimensional systems. These results were then 
  extended in \cite{ammar2009kalman}, where a necessary and sufficient 
  condition is given for constant coefficients but different diffusion coefficients $d_i$ 
  (the Laplace operator $\Delta$ can also be replaced by some general time-independent elliptic operator). 
  Moreover, in \cite{ammar2009generalization}, 
 some results are obtained in the case of time-dependent coefficients under  a sharp 
sufficient condition which is similar to the sharp \emph{Silverman-Meadows Theorem} 
in the finite-dimensional case. 

Concerning the case of space-varying coefficients, 
there is currently no general theory. In the case where the support of the coupling terms intersect 
the control domain, the most general result is proved in \cite{gonzalez2010controllability} 
for parabolic systems in cascade form with one control force (and possibly one order coupling terms). 
We also mention  \cite{MR2226005}, where a result of null controllability is proved in 
the case of a system of two equations with one control force, with an application to the controllability 
of a non-linear system of transport-diffusion equations.
In the case where the coupling regions do not intersect the control domain, 
only few results are known and in general there are some technical and geometrical restrictions 
(see for example \cite{alabau2013}, \cite{MR3039207} or \cite{MR2783322}). 
These restrictions come from the use of the \emph{transmutation method} that requires
 a controllability result on some related hyperbolic system. 
Let us also mention \cite{CherifMinimalTimeDisjoint}, where the authors consider a system of two equations in one space 
dimension and obtain a minimal time for null controllability, when the supports of the coupling terms 
do not intersect the control domain.

Concerning the case of first order coupling terms, there are also only few results.
The first one is \cite{guerrerosyst22}, where the author  studies notably
the case of $N+1$ coupled heat equations with $N$ control forces 
(we recall that $N$ is the dimension of $\Omega$) 
and  obtains the null controllability of System \eqref{system primmal} at any time 
when  the following estimate holds:
\begin{equation}\label{cond sergio}
 \|u\|_{H^1(\Omega)}\leqslant C\|(g_{21}\cdot\nabla +a_{21})^*u\|_{L^2(\Omega)^N},
\end{equation}
for all $u\in H^1_0(\Omega)$. Let us emphasize that inequality \eqref{cond sergio} 
is very restrictive, because it notably implies that $g_{21}$ has to be non-zero on each of its components (due to the $H^1-$norm appearing in the left-hand side).
Moreover, in the case of two equations, the result given in \cite{guerrerosyst22} 
is true only in the dimension one. 
Another case is the null controllability at any time $T>0$ of $m$ equations with one control force, which is studied 
in  \cite{gonzalez2010controllability}, under many assumptions: the coupling matrix $G$ has 
to be upper triangular (except on the controlled equation) and many coefficients of $A$ 
have to be non identically equal to zero on $\omega$, notably, in the $2\times2$ case, we should have
\begin{equation}\label{cond casade}\begin{array}{c}
g_{21}\equiv0 \mathrm{~in~} q_T \\
\mathrm{and}\\
(a_{21}>a_0 \mathrm{~in~} q_T \mathrm{~~~or~~~}a_{21}<-a_0 \mathrm{~in~} q_T), 
\end{array}\end{equation}
for a constant $a_0>0$.
The last result concerning first order coupling terms is the recent work \cite{benabdallah2014}, 
where the case of $2\times2$ and $3\times 3$ systems with one control force is studied under some technical assumptions. Notably, in the   $2\times 2$ case, the authors assume that 
 \begin{equation}\label{cond bord}
 \left\{\begin{array}{l}
\mathrm{there~exists~an~nonempty~open~subset}~\gamma~\mathrm{of}~ \partial\omega\cap\partial\Omega,\\
 \exists x_0\in\gamma \mathrm{~s.t.~}g_{21}(t,x_0)\cdot \nu(x_0)\neq0 \mathrm{~~for~all~}t\in[0,T],
\end{array}\right.
 \end{equation}
 where $\nu$ represents the exterior normal unit vector  to the boundary $\partial\Omega$. 
Under these technical restrictions on the control domain and the coupling terms, 
System \eqref{system primmal} is null controllable at any  time $T>0$.

Here we detail how our results differ from the existing ones: 
\begin{enumerate}
\item In the case of constant coefficients, we are able to obtain a necessary 
and sufficient condition in the case of $m$ equations, $m-1$ controls 
and coupling terms of order $0$ or $1$, 
which is the main new result.  Moreover, 
the diffusion coefficients can be different, 
we are able to treat the case of  as many equations as wanted 
 and we use an inequality  similar to 
Condition \eqref{cond sergio} but with the  $L^2-$norm in the left-hand side (see Lemma \ref{poincare}). 
To finish, we do not need  the control domain to extend up to the boundary as in Condition \eqref{cond bord}. 
The main restriction is that all the coefficients of System \eqref{system primmal} must be constant.
\item In the case  $N=1$ and $m=2$, we are able to obtain the controllability in arbitrary 
small time under some generic condition which is  purely technical. 
In Theorem 2, the geometric condition \eqref{cond bord} is not necessary, which is satisfying.
\end{enumerate}
\subsection{Related open problems and perspectives}
Let us describe briefly some related open problems and the difficulties that prevent us to go further in the present paper:
\begin{itemize}
\item As explained in Remark $1$, we believe that the condition involving the determinant of $H$ is purely technical, which means that either the first line 
of \eqref{cond1th2} or \eqref{g21 non nul} is sufficient as soon as we restrict 
to the class of coupling terms that intersect the control region. This conjecture includes the particular cases mentioned in item (e) of Remark $1$ (for instance the condition appearing in the case of time-dependent coefficients should not be necessary). However, we were not able to treat the general case because it is \textbf{crucial} in the proof of Theorem $1$ that the coefficients are \textbf{constant}, at least at two levels:
\begin{enumerate}
\item Lemma \ref{poincare} is not true anymore for variable coefficients (maybe in the case of time-varying coefficients only this can be adapted though).
\item During the proof of the Carleman estimate \eqref{ine obs 1}, it is crucial that the operator $\mathcal N$ introduced in \eqref{def N} has constant coefficients, so that it commutes with all the terms of the equation. This enables us to obtain an equation like \eqref{dual 2x2 psi chap B}, hence we apply directly Lemma \eqref{carleman 2X2 constant}, which would not be possible in the case of 
non-constant coefficients (there will be some remaining terms that cannot be absorbed).
\end{enumerate}
\item Another important limitation of our study is that we are restricted to the case of $m-1$ controls. One natural extension would then be to try to remove more controls. However, even in the case of $m$ equations and $m-2$ controls ($m\geqslant 3$) and constant coefficients, the situation is much more intricate because we do not have any idea on what would be a ``natural condition'' similar to the one of Theorem $1$ that will be enough to ensure the null controllability. Moreover, using the fictitious control method is not totally straightforward, because it is likely that the operator $\mathcal N$ needed in \eqref{def N} has to be of order $2$ in space, which would lead to important technical issues since it would require to prove a Carleman estimate similar to \eqref{ine obs 1} with a local term involving derivatives of order $2$.

\item 
To prove the \textit{algebraic solvability} 
 of system \eqref{syst21} (see \eqref{LMB} for the definition of this notion), 
we introduce a differential operator $\mathcal Q$ in \eqref{defQ2} which is a key issue of the proof
of the second item in   Theorem \ref{th 2}. 
The matrix $H$ given in \eqref{defH} 
is then deduced from $\mathcal Q$. The operator $ \mathcal Q$ has been found with the help of formal calculations. 
The authors have not been able to find  simpler expressions.  
It could be interesting to improve it and notably to investigate if looking at higher derivatives would give some different (and simpler) conditions.
\end{itemize}

\subsection{Strategy}\label{section strategy}

 The method described in this section is sometimes called  \textit{fictitious control method} 
and has already been used for instance  in \cite{gonzalezperez2006}, \cite{coronlissy2014} and \cite{ACO}. One important limitation of this method is that it will never 
be useful to treat the case where the support of the coupling terms do not intersect the control region, 
because, in what we call the algebraic resolution, we have to work locally on the control region.

Roughly, the method is the following: we first control the equations with $m$ controls (one on each equation) and we try to eliminate the control on the last equation thanks to algebraic manipulations.
Let us be more precise and decompose the problem into two different steps:

\vspace*{0.5cm}

\textbf{Analytic problem:} \\
Find a solution  $(z,v)$ in an appropriate space 
to the control problem by $m$ controls which are regular enough and are in the range of a differential operator. 
More precisely, solve
\begin{equation}\label{solution deux controle}
 \left\{\begin{array}{ll}
       \partial_tz=\Div(D\nabla z)
       +G\cdot\nabla  z+Az+\mathcal{N}(\mathds{1}_{\widetilde\omega}v)&\mathrm{in}~ Q_T,\\
                 z=0&\mathrm{on}~\Sigma_T,\\
        z(0,\cdot)=y^0,~z(T,\cdot)=0&\mathrm{in}~\Omega,
        \end{array}
\right.
\end{equation}
where $\mathcal{N}$ is some differential operator to be chosen later and $\widetilde\omega$ is strongly included in $\omega$.
Solving Problem \eqref{solution deux controle} is easier than solving the null controllability 
at time $T$ of System \eqref{system primmal}, 
because we control System \eqref{solution deux controle}  with a control on each equation.  
The important points (and somehow different from the usual methods) are that:
\begin{enumerate}
\item The control has to be of a special form 
(it has to be in the range of a differential operator $\mathcal{N}$),
\item The control 
has to be regular enough, so that it can be differentiated a certain amount 
of times with respect to the space and/or time variables (see the next section about the algebraic resolution).
\end{enumerate}
If we look for a control $v$ in the weighted Sobolev space $L^2(Q_T,\rho^{-1/2})^m$ for some weight $\rho$, 
it is well known (see, e.g. \cite[Th. 2.44, p. 56–57]{coron2009control}) 
that the null controllability at time $T$ of System \eqref{solution deux controle} 
  is equivalent to the following \textit{observability inequality}:
 \begin{equation}\label{ine obs strat}
 \begin{array}{c}
 \displaystyle\int_{\Omega} 
|\psi(0,x)|^2dx
\leqslant C_{obs} 
\displaystyle\iint_{q_T} \rho
|\mathcal{N}^*\psi(t,x)|^2\,dxdt,
\end{array} \end{equation}
where $\psi$ is the solution to the dual system
\begin{equation*}
 \left\{\begin{array}{ll}
-\partial_t\psi=\Div(D\nabla  \psi)-G^*\cdot\nabla  \psi+A^*\psi&\mathrm{in}~ Q_T,\\
       \psi=0&\mathrm{on}~\Sigma_T,\\
       \psi(T,\cdot)=\psi^0&\mathrm{in}~\Omega.
        \end{array}
\right.
\end{equation*}
Inequalities like \eqref{ine obs strat} can be proved thanks 
to some appropriate \emph{Carleman estimates}. The weight $\rho$ can be chosen to be 
exponentially decreasing at times $t=0$ and $t=T$, which will be useful later. 
In fact, we will have to adapt  the usual HUM duality method to 
ensure that one can find such controls.

\vspace*{0.5cm}

\textbf{Algebraic problem:} \\
For $f:=\mathds{1}_{\omega}v$,   
 find a pair $(\widehat z,\widehat v)$ (where $\widehat v$ now acts only on the first $m-1$ equations) in an appropriate space satisfying the following control problem:
\begin{equation}\label{probleme ramene a tout l espace}
 \left\{\begin{array}{ll}
       \partial_t\widehat{z}=\Div(D\nabla \widehat{z})+G\cdot\nabla \widehat{z}+A\widehat{z}+B\widehat{v}+\mathcal{N}f&\mathrm{in}~ Q_T,\\
       \widehat{z}=0&\mathrm{on}~\Sigma_T,\\
       \widehat{z}(0,\cdot)=\widehat{z}(T,\cdot)=0&\mathrm{in}~\Omega,
        \end{array}
\right.
\end{equation}
and such that the spatial support of $\widehat v$
is strongly included in $\omega$.  
We will solve this problem using the notion of \emph{algebraic solvability} 
of differential systems, which is based on ideas coming from \cite[Section 2.3.8]{Gromovbook} 
and was already widely used in \cite{coronlissy2014} and \cite{ACO}. 
The idea is to write System \eqref{probleme ramene a tout l espace} as an \emph{underdetermined} 
 system in the variables $\widehat{z}$ and $\widehat{v}$ and to see $\mathcal{N}f$ as a source term, 
 so that we can write Problem \eqref{probleme ramene a tout l espace}  under the abstract form
\begin{equation}\label{def L Strategy 2}
\mathcal{L}(\widehat{z},\widehat{v})=\mathcal{N}f,
\end{equation}
where 
\begin{equation}\label{def L Strategy}
 \mathcal{L}(\widehat{z},\widehat{v}):=\begin{array}{ll}
       \partial_t\widehat{z}-\Div(D\nabla \widehat{z})-G\cdot\nabla \widehat{z}-A\widehat{z}-B\widehat{v}.
        \end{array}
\end{equation}
The goal will be then to find a partial differential operator $\mathcal{M}$ satisfying
\begin{gather}
\label{LMB}
\mathcal L \circ \mathcal M=\mathcal N.
\end{gather}
When \eqref{LMB} is satisfied, we say that  System \eqref{probleme ramene a tout l espace}
is \emph{algebraically solvable}. 
This exactly means that one can find a solution $(\widehat{z},\widehat{v})$ to System  
\eqref{probleme ramene a tout l espace} which can be written 
as a linear combination of some derivatives of the source term $\mathcal N f$. 
The main advantage of this method is that one can only work \emph{locally} on 
$q_T$, because the solution depends locally on the source term and then has 
the same support as the source term (to obtain a solution which is defined everywhere on $Q_T$, 
one just extends it by $0$). This part will be explained in more details in Sections $2.1$ and $3.1$.

\vspace*{0.5cm}

\textbf{Conclusion:}\\
If we can solve the analytic and algebraic problems, 
then it is easy to check that   $(y,u):=(z-\widehat{z},-\widehat{v})$ will be a solution to 
System (\ref{system primmal}) 
in an appropriate space and will satisfy $y(T)\equiv0$ in $\Omega$ 
(for more explanations, see \cite[Prop. 1]{coronlissy2014} or Sections 2.4 and 3.3 of the present paper).

\section{Proof of Theorem \ref{th 1}}

Let us remind that in this case, $D$, $G$ and $A$ are constant. 
In Section 2.1, 2.2 and 2.3, we will always consider some $i_0\in\{1,...,m-1\}$ such that 
\begin{equation}\label{cond ine}
 g_{mi_0}\neq0 \mathrm{~or~}a_{mi_0}\neq0. 
\end{equation}
We will follow the strategy described in Section \ref{section strategy}, and we first begin with finding some appropriate operator $\mathcal N$.

\subsection{Algebraic resolution}\label{sec alg resol const}

 We will here explain how to choose the differential operator   $\mathcal{N}$ used in the next section. 
We will assume from now on that all differential operators of this section are 
defined in $\mathcal{C}^{\infty}(Q_T)^m$. The appropriate spaces will be specified in Section  \ref{sec proof th1}. 
 We consider $\mathcal{N}$ as the operator defined for all $f:=(f_1,\ldots,f_m)$ by
\begin{equation}\label{def N}
\mathcal{N}(f):= \left(\begin{array}{c}
\mathcal{N}_1f\\
\mathcal{N}_2f\\
\ldots\\
\mathcal{N}_mf
        \end{array}\right):=\left(\begin{array}{c}(-g_{mi_0}\cdot\nabla-a_{mi_0})f_1\\(-g_{mi_0}\cdot\nabla-a_{mi_0})f_2\\
       \ldots\\(-g_{mi_0}\cdot\nabla-a_{mi_0})f_m\end{array}\right).\end{equation}
 Let us recall that the definition of $\mathcal{L}$ is given in \eqref{def L Strategy}.
 
 As explained in Section \ref{section strategy},
we want find a differential operator   $\mathcal{M}$ such that 
\begin{equation}\label{LM=N}
 \mathcal{L}\circ\mathcal{M}=\mathcal{N}.
\end{equation}
 We have the following proposition:
 \begin{prop}\label{alg1} Let $\mathcal N$ be defined as in \eqref{def N}. Then there exists a differential operator 
        $\mathcal{M}$ of order $1$ in time and $2$ in space, with constant coefficients, such that \eqref{LM=N} is verified. 
\end{prop}

\noindent\textbf{Proof of Proposition \ref{alg1}.}
We can remark that equality \eqref{LM=N} is equivalent to
\begin{equation}\label{ M*L*=N*}
\mathcal{M}^*\circ \mathcal{L}^*=\mathcal{N}^*.
\end{equation}
The adjoint $\mathcal{L}^*$ of the operator $\mathcal{L}$ is given for all $\varphi\in \mathcal{C}^{\infty}(Q_T)^m$ by
\begin{equation}\label{def L*}
\begin{array}{rcl}
 \mathcal{L}^*\varphi&:=& \left(\begin{array}{c}
\mathcal{L}^*_1\varphi\\
\ldots\\
\mathcal{L}^*_{2m-1}\varphi
        \end{array}\right)
= \left(\begin{array}{c}
-\partial_t\varphi_1-\Div(d_1\nabla \varphi_1)+\sum_{j=1}^{m}\{g_{j1}\cdot\nabla \varphi_j-a_{j1}\varphi_j\}\\
\ldots\\
-\partial_t\varphi_m-\Div(d_m\nabla \varphi_m)+\sum_{j=1}^{m}\{g_{jm}\cdot\nabla \varphi_j-a_{jm}\varphi_j\}
\\
\varphi_1\\
\ldots\\
\varphi_{m-1}

        \end{array}\right),
\end{array}\end{equation}
the $m-1$ last lines coming from the particular form of our control operator $B$.
Now we apply $g_{mi_0}\cdot\nabla-a_{mi_0}$ to the $(m+i)^{th}$ line for $i\in \{1,...,m-1\}$ and  we add
$(\partial_t+\Div(d_{i_0}\nabla))\mathcal{L}^*_{m+i_0}\varphi
+\sum_{j=1}^{m-1}(-g_{ji_0}\cdot\nabla+a_{ji_0})\mathcal{L}^*_{m+j}\varphi$ to the $(i_0)^{th}$ line. 
Hence, remarking that $\mathcal{L}^*_{m+i}\varphi=\varphi_i$ ($i\in \{1,...,m-1\}$), we obtain
$$
 \left(\begin{array}{c}
(g_{mi_0}\cdot\nabla-a_{mi_0})\mathcal{L}^*_{m+1}\varphi\\
\ldots\\ (g_{mi_0}\cdot\nabla-a_{mi_0})\mathcal{L}^*_{2m-1}\varphi
\\
\mathcal{L}^*_{i_0}\varphi
+(\partial_t+\Div(d_{i_0}\nabla))\mathcal{L}^*_{m+i_0}\varphi
+\sum_{j=1}^{m-1}(-g_{ji_0}\cdot\nabla+a_{ji_0})\mathcal{L}^*_{m+j}\varphi
        \end{array}\right)  =
\mathcal{N}^*\varphi.
$$
Thus if we define $\mathcal{M}^*$ for $\psi:=(\psi_1,...,\psi_{2m-1})\in \mathcal{C}^{\infty}(Q_T)^{2m-1}$ by
\begin{equation}\label{def M*}
\begin{array}{l}
 \mathcal{M}^* \left(\begin{array}{c}
\psi_1\\\ldots\\\psi_{2m-1}
        \end{array}\right)
        :=\left(\begin{array}{c}
(g_{mi_0}\nabla-a_{mi_0})\psi_{m+1}\\\ldots\\
(g_{mi_0}\nabla-a_{mi_0})\psi_{2m-1}\\
\psi_{i_0}
+(\partial_t+\Div(d_{i_0}\nabla))\psi_{m+i_0}+\sum_{j=1}^{m-1}(-g_{ji_0}\cdot\nabla+a_{ji_0})\psi_{m+j}
        \end{array}\right),
                 \end{array}
        \end{equation}
then equality \eqref{ M*L*=N*} is satisfied and hence equality \eqref{LM=N} also. Moreover, the coefficients of $\mathcal M^*$ are constant, hence it is also the case for the coefficients of $\mathcal M$.
\cqfd

\begin{rem}Looking carefully at the proof of Proposition \ref{alg1}, we remark that one could also have constructed some differential operator $\widetilde {\mathcal M}$ such that  $\mathcal{L}\circ\widetilde {\mathcal M}=\widetilde {\mathcal N}$ with $\widetilde {\mathcal N}_1=\ldots=\widetilde {\mathcal N}_{m-1}=Id$ and $\widetilde {\mathcal N}_m=-g_{mi_0}.\nabla -a_{mi_0}$. However, it is more convenient for the proof of \eqref{ine obs 1} to work with a differential operator of same order on each line of $\mathcal N$ and that is the reason why we copied $-g_{mi_0}.\nabla -a_{mi_0}$ on each line.
\end{rem}

\subsection{An appropriate Carleman estimate}\label{sec ine obs}

 Let us consider the following dual system associated to System \eqref{solution deux controle}
\begin{equation}\label{primal 2x2 constant dirch}
 \left\{\begin{array}{ll}
-\partial_t\psi=\Div(D\nabla \psi)-G^*\cdot\nabla \psi+A^*\psi&\mbox{ in }~ Q_T,\\
 \psi=0&\mbox{ on }~\Sigma_T,\\
   \psi(T,\cdot)=\psi^0&\mbox{ in }~\Omega.
        \end{array}
\right.
\end{equation}

The two main results of this section are Propositions \ref{prop ine obs 1} and \ref{prop ine obs 2}, 
which are respectively some Carleman estimate and observability inequality. 
The particularity of theses inequalities is that the observation will not be directly 
the $L^2-$norm of the solution $\psi$ to System \eqref{primal 2x2 constant dirch} on the subset $\omega$, 
but it will be the $L^2-$norm of some linear combination of $\psi$ and its derivatives of first order on the subset $\omega$. 
This particular form
will be used in the next section to construct a solution to the analytic 
control Problem \eqref{solution deux controle}. 

Let $\omega_0$, $\omega_1$ and $\omega_2$ be three nonempty open subsets  included in $\omega$ satisfying
\begin{equation*}
\overline{\omega_2}\subset\omega_1,~
\overline{\omega_1}\subset\omega_0\mathrm{~and~}
\overline{\omega_0}\subset\omega.
\end{equation*}
Before stating the Carleman estimate, let us introduce some notations. For $s,\lambda>0$, let us define
\begin{equation}\label{defI}
 I(s,\lambda;u):=
 s^3\lambda^4\displaystyle\iint_{Q_T} e^{-2s\alpha}\xi^3|u|^2dxdt
 + s\lambda^2\displaystyle\iint_{Q_T} e^{-2s\alpha}\xi|\nabla u|^2dxdt,
\end{equation}
where 
\begin{equation}\label{defax}
 \alpha(t,x):=\dfrac{\mathrm{exp}(12\lambda\|\eta^0\|_{\infty})
 -\mathrm{exp}[\lambda(10\|\eta^0 \|_{\infty}+\eta^0(x) )]}{t^5(T-t)^5}
 ~~~\mathrm{and}~~~
 \xi(t,x):=\dfrac{\mathrm{exp}[\lambda(10\|\eta^0 \|_{\infty}+\eta^0(x) )]}{t^5(T-t)^5}.
\end{equation}
Here, $\eta^0 \in\mathcal{C}^2(\overline{\Omega})$ is a function satisfying
\begin{equation*}
 |\nabla\eta^0 |\geqslant \kappa\mathrm{~in~}\Omega\backslash\omega_2,~~~
  \eta^0>0\mathrm{~in~}\Omega ~~~\mathrm{and}~~~
   \eta^0=0\mathrm{~on~}\partial\Omega,
\end{equation*}
with $\kappa>0$. The proof of the existence of such a function $\eta^0 $  
can be found  in \cite[Lemma 1.1, Chap. 1]{fursikov1996controllability} 
(see also \cite[Lemma 2.68, Chap. 2]{coron2009control}). We will use the two notations 
\begin{equation}\label{defaxs}
 \alpha^*(t):=\max\limits_{x\in\overline{\Omega}}\alpha(t,x)
 \mathrm{~~~and~~~}\xi^*(t):=\min\limits_{x\in\overline{\Omega}}\xi(t,x),
\end{equation}
for all $t\in (0,T)$.

\subsubsection{Some auxiliary results}

Let us now give some useful auxiliary results that we will need in our proofs. 
The first one is  a Carleman estimate which holds for solutions of the heat equation with 
non-homogeneous Neumann boundary conditions:

\begin{Lemme}\label{th carl neum} 
 Let us assume that $d>0$, $u^0\in L^2(\Omega)$, $f_1\in L^2(Q_T)$ and $f_2\in L^2(\Sigma_T)$. 
 Then there exists a constant $C:=C(\Omega,\omega_2)>0$ such that the solution to the system 
 \begin{equation*}\left\{\begin{array}{ll}
  -\partial_tu-\Div(d\nabla u) =f_1&\mathrm{in}~Q_T,\\
  \frac{\partial u}{\partial n} =f_2&\mathrm{on}~\Sigma_T,\\
  u(T,\cdot)=u^0&\mathrm{in~}\Omega,\end{array}\right.
 \end{equation*}
satisfies
\begin{equation*}\begin{array}{r}
 I(s,\lambda;u)\leqslant C\left(
 s^3\lambda^4\displaystyle\iint_{(0,T)\times\omega_2} e^{-2s\alpha}\xi^3|u|^2dxdt
 +\displaystyle\iint_{Q_T} e^{-2s\alpha}|f_1|^2dxdt
 ~~~~~~~~~~~~~~~~~~~~~\right.\\\left.
 +s\lambda\displaystyle\iint_{\Sigma_T}e^{-2s\alpha^*}\xi^*|f_2|^2 d\sigma dt\right),
\end{array}\end{equation*}
for all $\lambda\geqslant C$ and $s\geqslant C(T^5+T^{10})$.
\end{Lemme}
The proof of Lemma \ref{th carl neum} can essentially be found in \cite{GuerreroFourier}. 
In fact, in this article, 
the weights are a little bit different ($t(T-t)$ instead of $t^5(T-t)^5$),  
but the proof just needs to be slightly adapted to obtain the present result.

From Lemma \ref{th carl neum}, one can deduce the following result:

\begin{Lemme}\label{carleman 2X2 constant}
 Let  $h\in L^2(\Sigma_T)^m$. 
Then there exists a constant $C:=C(\Omega,\omega_2)>0$ such that 
for every $\varphi^0\in L^2(\Omega)^m$, 
the solution $\varphi$ 
 to the system 
\begin{equation}\label{primal 2x2 constant fourier}
 \left\{\begin{array}{ll}
-\partial_t\varphi=\Div(D\nabla \varphi)
-G^*\cdot\nabla  \varphi+A^*\varphi&\mathrm{in}~ Q_T,\\
       \frac{\partial\varphi}{\partial n}=h&\mathrm{on}~\Sigma_T,\\
       \varphi(T,\cdot)=\varphi^0&\mathrm{in}~\Omega
        \end{array}
\right.
\end{equation}
satisfies
\begin{equation*}\begin{array}{l}
 I(s,\lambda;\varphi)
 \leqslant C\left(
 s^3\lambda^4\displaystyle\iint_{(0,T)\times\omega_2} e^{-2s\alpha}\xi^3|\varphi|^2dxdt
 +s\lambda\displaystyle\iint_{\Sigma_T}e^{-2s\alpha^*}\xi^*\left|h\right|^2d\sigma dt \right),
\end{array}\end{equation*}
for every $\lambda\geqslant C$ and $s\geqslant s_0=C(T^5+T^{10})$.
\end{Lemme}


The proof of Lemma \ref{carleman 2X2 constant} is standard and is left to the reader (one just have to apply Lemma \ref{th carl neum} separately to all equations of System \eqref{primal 2x2 constant fourier}, sum of all the Carleman estimates and absorb the remaining lower-order terms thanks to the left-hand side).

\vspace{0,2cm}

In this section, we will use also the following estimate.

\begin{Lemme}\label{poincare poids}
Let $r\in \mathbb{R}$. Then there exists $C:=C(r,\omega_2,\Omega)>0$  
such that, for every $T>0$ and every $u\in L^2(0,T;H^1(\Omega))$,
\begin{equation*}\begin{array}{r}
 s^{r+2}\lambda^{r+2}\displaystyle\iint_{Q_T} e^{-2s\alpha}\xi^{r+2}|u|^2dxdt
 \leqslant C\left(s^r\lambda^r
  \displaystyle\iint_{Q_T} e^{-2s\alpha}\xi^{r}|\nabla u|^2dxdt
  \right.~~~~~~~~~~~~~~~~~~~~~~~\\\left.
  + s^{r+2}\lambda^{r+2}\displaystyle\iint_{(0,T)\times\omega_2} e^{-2s\alpha}\xi^{r+2}|u|^2dxdt\right),
\end{array}\end{equation*}
for every $\lambda\geqslant C$ and $s\geqslant C(T^5+T^{10})$.
\end{Lemme}

The proof of this lemma can be found for example in \cite[Lemma 3]{coronguerrero2009}. 
Our next Lemma is some Poincar\'e-type inequality involving  the differential operator $\mathcal{N}^*$.

\begin{Lemme}\label{poincare}
There exists a constant $C:=C(\Omega)>0$ 
such that  for every $u\in H_0^1(\Omega)$, the following estimate holds:
\begin{equation}\label{ine poincare ame}
 \displaystyle\int_{\Omega}u^2
 \leqslant C\displaystyle\int_{\Omega}|\mathcal{N}^*u|^2,
\end{equation}
where $\mathcal{N}^*:=g_{mi_0}\cdot\nabla-a_{mi_0}$.
\end{Lemme}

Lemma \ref{poincare} is obvious if $g_{mi_0}=0$. If $a_{mi_0}=0$, this is exactly the usual Poincar\'e inequality. The case  $a_{mi_0}\not =0$ and $g_{mi_0}\not = 0$ can be reduced to the previous case by considering $$u(x)\exp\left(-\frac{a_{mi_0}}{||g_{mi_0}||^2}(g_{mi_0}\cdot x)\right).$$

\vspace{0,2cm}

In order to deal with more regular solutions, one needs the following lemma:

\begin{Lemme} \label{lemme regul}
Let  $z_0\in H^1_0(\Omega)^m$ and $f\in  L^2(Q_T)^m$. 
Let us denote by $\mathcal{R}:=-\Div(D\nabla\cdot)-G\cdot\nabla-A$ and 
  consider $z$ the solution in $W^{2,1}_2(Q_T)^m$ to the system
    \begin{equation}\label{systeme regul}
 \left\{\begin{array}{ll}
\partial_tz=\Div(D\nabla z)+G\cdot\nabla z+Az
+f&\mathrm{in}~ Q_T,\\
       z=0&\mathrm{on}~\Sigma_T,\\
       z(0,\cdot)=z_0&\mathrm{in}~\Omega.
        \end{array}
\right.
\end{equation}
Let $d\in\mathbb N$. Let us assume that 
$z_0\in H^{2d+1}(\Omega)^m$, $f\in L^2(0,T;H^{2d}(\Omega)^m)\cap H^{d}(0,T;L^{2}(\Omega)^m)$ 
and satisfy the following compatibility conditions:
\begin{equation*}
 \left\{\begin{array}{l}
      g_0:=z_0\in H^1_0(\Omega)^m,  \\
      g_1:= f(0,\cdot)-\mathcal{R}g_0\in H^1_0(\Omega)^m, \\
      \ldots\\
      g_{d}:= \partial^{d-1}_t f(0,\cdot)-\mathcal{R}g_{d-1}\in H^1_0(\Omega)^m.
        \end{array}
\right.
\end{equation*}
 Then $z\in L^2(0,T;H^{2d+2}(\Omega)^m)\cap H^{d+1}(0,T;L^{2}(\Omega)^m)$ and we have the estimate
 \begin{equation}\label{estim regul 1}
  \|z\|_{ L^2(0,T;H^{2d+2}(\Omega)^m)\cap H^{d+1}(0,T;L^{2}(\Omega)^m)}
  \leqslant C(\|f\|_{L^2(0,T;H^{2d}(\Omega)^m)\cap H^d(0,T;L^{2}(\Omega)^m)}
  +\|z_0\|_{H^{2d+1}(\Omega)^m}).
 \end{equation}

\end{Lemme}
It is a classical result that can be easily deduced for example from \cite[Th. 6, p. 365]{MR2597943}.

\subsubsection{Carleman inequality}
We are now able to prove the following inequality:

\begin{prop}\label{prop ine obs 1}
 There exists a   constant $C:=C(\omega_0,\Omega)>0$ such that for every 
  $\psi^0\in L^2(\Omega)^m$, 
 the corresponding solution $\psi$ to System \eqref{primal 2x2 constant dirch} satisfies
 \begin{equation}\label{ine obs 1}
 \begin{array}{c}
\displaystyle\iint_{Q_T} e^{-2s\alpha}\{
s^7\lambda^8\xi^{7}|\mathcal{N}^*\psi|^2+s^5\lambda^6\xi^{5}|\nabla\mathcal{N}^*\psi|^2
+s^3\lambda^4\xi^{3}|\nabla\nabla\mathcal{N}^*\psi|^2+s\lambda^2\xi|\nabla\nabla\nabla\mathcal{N}^*\psi|^2   \}dxdt\\
\leqslant C s^7\lambda^8
\displaystyle\iint_{(0,T)\times\omega_0} e^{-2s\alpha}\xi^7
|\mathcal{N}^*\psi|^2dxdt,
\end{array} 
\end{equation}
for every $\lambda \geqslant C$ and 
$s\geqslant s_0= C(T^5+T^{10})$. 
\end{prop}
\begin{rem}
It may be quite surprising that one can put so much derivatives at the left-hand-side of equality \eqref{ine obs 1}, because the initial condition $\psi^0$ is only supposed to be $L^2$, 
hence $\psi$ is only assumed to be in $W(0,T)$ (see \eqref{defW}). However, because of the fact that the exponential weight $e^{-2s\alpha}$ is strong enough to 
absorb the singularity that only exists at initial time $t=0$, it is quite easy to prove that all the integrals appearing in the left-hand side of \eqref{ine obs 1} exist 
(this can notably be deduced for example from inequalities like \eqref{estim41}, \eqref{estim43} or \eqref{estim61}).
\end{rem}
\noindent\textbf{Proof of Proposition \ref{prop ine obs 1}.}

The proof is inspired by \cite{coronguerrero2009}. The main difference here is that we keep $\mathcal N^* \psi$ at the right-hand side, which complicates a little bit the proof.
 Let us denote by 
 \begin{equation}\label{def mcR}
\mathcal{R}:= -\Div(D\nabla)+G^*\cdot\nabla-A^*. 
 \end{equation}
 We can assume without loss of generality that 
\begin{equation*}
\psi^0\in H^5(\Omega)\mathrm{~~and~~}
\psi^0,\mathcal{R}\psi^0,\mathcal{R}^2\psi^0\in H^1_0(\Omega)^m
\end{equation*}
(The general case follows from a density argument). Thus, using Lemma \ref{lemme regul}, 
the solution $\psi$ to System \eqref{primal 2x2 constant dirch} is an element of 
$ L^2(0,T;H^{6}(\Omega)^m)\cap H^{3}(0,T;L^{2}(\Omega)^m)$. 
First of all, let us  apply the differential operator 
$$\nabla\nabla\mathcal{N}^*=\nabla\nabla(-a_{mi_0}+g_{mi_0}\cdot\nabla)$$ 
to System \eqref{primal 2x2 constant dirch} 
satisfied by $\psi$. Thus, if we call  $\phi:=(\phi_{ij})_{1\leqslant i,j\leqslant N}$ 
with $\phi_{ij}:=\partial_i\partial_j\mathcal{N}^*\psi$, then one observes that $\phi$ is a solution of the following system:
 \begin{equation}\label{dual 2x2 psi chap B}
 \left\{\begin{array}{ll}
-\partial_t\phi_{ij}=\Div(D\nabla \phi_{ij})-G^*\cdot\nabla\phi_{ij}+A^*\phi_{ij}&\mathrm{in}~ Q_T,\\
       \frac{\partial\phi_{ij}}{\partial n}
       =\frac{\partial(\partial_i\partial_j\mathcal{N}^*\psi_{ij})}{\partial n}&\mathrm{on}~\Sigma_T,\\
       \phi_{ij}(T,\cdot)=\partial_i\partial_j\mathcal{N}^*\psi_{ij}^0&\mathrm{in}~\Omega.
        \end{array}
\right.
\end{equation}
By applying Lemma \ref{carleman 2X2 constant} to $\phi$, 
we have
\begin{equation}\label{estim preuve lemme carl2}
\begin{array}{l}
 I(s,\lambda,\phi)
 \leqslant C\left(
 s^3\lambda^4\displaystyle\iint_{(0,T)\times\omega_2} e^{-2s\alpha}\xi^3|\phi|^2dxdt
 +s\lambda\displaystyle\iint_{\Sigma_T}e^{-2s\alpha^*}\xi^*
 \left|\frac{\partial (\nabla\nabla \mathcal{N}^*\psi)}{\partial n}\right|^2
 d\sigma dt \right),
\end{array}\end{equation}
for every $\lambda \geqslant C$ and $s\geqslant C(T^5+T^{10})$. 

The  proof will be divided into three steps :
\begin{itemize}
 \item In the first step, we will estimate the boundary term in the right-hand side 
 of inequality \eqref{estim preuve lemme carl2}
with some global interior term involving $\psi$ that will be absorbed later.
 \item In the second step, we will compare $ I(s,\lambda,\phi)$ with the left-hand side of inequality  
 \eqref{ine obs 1}.
 \item Finally, in the last step, we will estimate the local term of high order appearing 
 in  inequality \eqref{estim preuve lemme carl2}  thanks to some local terms that will be absorbed in the 
 left-hand side of  inequality  \eqref{estim preuve lemme carl2}  and also thanks to the local term of
the right-hand side of inequality \eqref{ine obs 1}.
\end{itemize}

\textbf{Step 1:} 
Let us consider a function $\theta\in\mathcal{C}^2(\overline{\Omega})$ 
such that 
\begin{equation*}
\frac{\partial\theta}{\partial n}=\theta=1  \mathrm{~on~} \partial\Omega.  
\end{equation*}
After an integration by parts of the boundary term, we obtain
\begin{equation*}\begin{array}{l}
s\lambda\displaystyle\int_{0}^Te^{-2s\alpha^*}\xi^*\displaystyle\int_{\partial\Omega}
 \left|\frac{\partial \phi}{\partial n}\right|^2d\sigma dt \\
=s\lambda\displaystyle\int_{0}^Te^{-2s\alpha^*}\xi^*\displaystyle\int_{\partial\Omega}
 \frac{\partial\phi}{\partial n}\nabla\phi\cdot\nabla\theta d\sigma dt \\
=s\lambda\displaystyle\int_{0}^Te^{-2s\alpha^*}\xi^*\displaystyle\int_{\Omega}
 \Delta\phi\nabla\phi\cdot\nabla\theta dx dt 
+s\lambda\displaystyle\int_{0}^Te^{-2s\alpha^*}\xi^*\displaystyle\int_{\Omega}
 \nabla(\nabla\theta\cdot\nabla\phi)\cdot\nabla\phi dxdt.
 \end{array}\end{equation*}
Using successively Cauchy-Schwarz inequality and Young's inequality, 
we deduce that
 \begin{equation}\label{estim preuve lemme carl3}
 \begin{array}{rcl}
s\lambda\displaystyle\int_{0}^Te^{-2s\alpha^*}\xi^*\displaystyle\int_{\partial\Omega}
 \left|\frac{\partial\phi}{\partial n}\right|^2d\sigma dt 
&\leqslant &C\lambda\displaystyle\int_{0}^Te^{-2s\alpha^*}
 \|(s\xi^*)^{4/5}\psi\|_{H^4(\Omega)^{m}}\|(s\xi^*)^{1/5}\psi\|_{H^5(\Omega)^{m}} dt\\
 &\leqslant &C\lambda\displaystyle\int_{0}^Te^{-2s\alpha^*}(s\xi^*)^{8/5}
 \|\psi\|_{H^4(\Omega)^{m}}^2 dt\\
&&~~~~~~~~~~~~~~~~~~~~~~
 +C\lambda\displaystyle\int_{0}^Te^{-2s\alpha^*}(s\xi^*)^{2/5}
 \|\psi\|_{H^5(\Omega)^{m}}^2 dt.
 \end{array}\end{equation}
 Let us introduce
 $\widehat{\psi}:=\rho\psi$ with $\rho\in \mathcal{C}^{\infty}([0,T])$ defined by 
$$\rho:= (s\xi^*)^{a}e^{-s\alpha^*},$$ for some $a\in\mathbb R$ to be chosen later. 

One remark that $\rho$ verifies $\partial^i_{t}\rho(0)=0$ for all $i\in\mathbb N$.  
 Then $\widehat{\psi}$ is solution to the system
   \begin{equation}\label{dual 2x2 psi rond}
 \left\{\begin{array}{ll}
-\partial_t\widehat{\psi}=\Div(D\nabla \widehat{\psi})-G^*\cdot\nabla\widehat{\psi}+A^*\widehat{\psi}
-\rho_t\psi&\mathrm{in}~ Q_T,\\
       \widehat{\psi}=0&\mathrm{on}~\Sigma_T,\\
       \widehat{\psi}(T,\cdot)=0&\mathrm{in}~\Omega.
        \end{array}
\right.
\end{equation}
 Lemma \ref{lemme regul} gives for $\widehat{\psi}$ the estimate   
\begin{equation}\label{boot1}\begin{array}{rcl}
  \|\widehat{\psi}\|_{L^2(0,T;H^{2d+2}(\Omega)^m)\cap H^{d+1}(0,T;L^{2}(\Omega)^m)}
  \leqslant C\|\rho_t\psi\|_{ L^2(0,T;H^{2d}(\Omega)^m)\cap H^{d}(0,T;L^{2}(\Omega)^m)},
 \end{array}\end{equation}
 for $d\in \{0,1,2\}$. 
 Using the definitions of $\xi^*$ and $\alpha^*$ given in \eqref{defaxs} and the particular form of $\rho$ chosen, we have 
 
 \begin{equation}\label{rhot}|\partial_t\rho|\leqslant CT(s\xi^*)^{a+6/5}e^{-s\alpha^*},\end{equation}

 \begin{equation}\label{rhott}|\partial_{tt}\rho|\leqslant CT^2(s\xi^*)^{a+12/5}e^{-s\alpha^*}\end{equation}
and
 \begin{equation}\label{rhottt}|\partial_{ttt}\rho|\leqslant CT^3(s\xi^*)^{a+18/5}e^{-s\alpha^*}.\end{equation}
 Using  inequality  \eqref{boot1} with $\rho:=e^{-s\alpha^*}(s\xi^*)^{4/5}$ and $d=1$, we obtain
  \begin{equation}\label{estim42}\begin{array}{c}
\displaystyle\int_{0}^Te^{-2s\alpha^*}(s\xi^*)^{8/5}
 \|\psi\|_{H^4(\Omega)^{m}}^2 dt
 \\\leqslant C\left(\displaystyle\int_{0}^T
 \|\partial_t(e^{-2s\alpha^*}(s\xi^*)^{4/5})\psi\|_{H^2(\Omega)^{m}}^2 dt+\displaystyle\int_{0}^T
 \|\partial_{t}(\partial_{t}(e^{-s\alpha^*}(s\xi^*)^{4/5})\psi)\|_{L^2(\Omega)^{m}}^2 dt\right).
 \end{array}\end{equation}
 Applying now inequality \eqref{boot1} with $\rho:=\partial_t(e^{-s\alpha^*}(s\xi^*)^{4/5})$ and $d=0$, we get
 \begin{equation}\label{estim41}\begin{array}{rcl}
\displaystyle\int_{0}^T
 \|\partial_t(e^{-s\alpha^*}(s\xi^*)^{4/5})\psi\|_{H^2(\Omega)^{m}}^2 dt+\displaystyle\int_{0}^T
 \|\partial_{t}(\partial_{t}(e^{-s\alpha^*}(s\xi^*)^{4/5})\psi)\|_{L^2(\Omega)^{m}}^2 dt~~~~~~~~~~~
 \\
 \leqslant C\displaystyle\int_{0}^T
 \|\partial_{tt}(e^{-s\alpha^*}(s\xi^*)^{4/5})\psi\|_{L^2(\Omega)^{m}}^2 dt.
 \end{array}\end{equation}
 Using \eqref{rhott} with $a=4/5$  together with \eqref{estim41} and \eqref{estim42}, we deduce
   \begin{equation}\label{estim43}\begin{array}{rcl}
\displaystyle\int_{0}^Te^{-2s\alpha^*}(s\xi^*)^{8/5}
 \|\psi\|_{H^4(\Omega)^{m}}^2 dt
 \leqslant CT^2\displaystyle\int_{0}^Te^{-2s\alpha^*}(s\xi^*)^{32/5}
 \|\psi\|_{L^2(\Omega)^{m}}^2 dt.
 \end{array}\end{equation}
 Using exactly the same proof and taking into account \eqref{rhottt}, one can also prove that 
    \begin{equation}\label{estim61}\begin{array}{rcl}
\displaystyle\int_{0}^Te^{-2s\alpha^*}(s\xi^*)^{-4/5}
 \|\psi\|_{H^6(\Omega)^{m}}^2 dt
 \leqslant CT^3\displaystyle\int_{0}^Te^{-2s\alpha^*}(s\xi^*)^{32/5}
 \|\psi\|_{L^2(\Omega)^{m}}^2 dt.
 \end{array}\end{equation}
Using  inequalities \eqref{estim43} and \eqref{estim61}, the interpolation inequality 
$$\|u\|_{H^5(\Omega)^{m}}\leqslant C \|u\|_{H^4(\Omega)^{m}}^{1/2}
\|u\|_{H^6(\Omega)^{m}}^{1/2}\mathrm{~for~every~}u\in H^6(\Omega)^{m},$$ and the Cauchy-Schwarz inequality, we deduce that
  \begin{equation}\label{estim preuve lemme carl5}\begin{array}{rcl}
\displaystyle\int_{0}^Te^{-2s\alpha^*}(s\xi^*)^{2/5}
 \|\psi\|_{H^5(\Omega)^{m}}^2 dt
  &\leqslant & C\displaystyle\int_{0}^T
 \|e^{-s\alpha^*}(s\xi^*)^{-2/5}\psi\|_{H^6(\Omega)^{m}}\|e^{-s\alpha^*}(s\xi^*)^{4/5}\psi\|_{H^4(\Omega)^{m}} dt\\
   &\leqslant & CT^{5/2}\displaystyle\int_{0}^Te^{-2s\alpha^*}(s\xi^*)^{32/5}\|\psi\|_{L^2(\Omega)^{m}}^2 dt.
\end{array}  \end{equation}
Thus inequalities \eqref{estim preuve lemme carl2}, \eqref{estim preuve lemme carl3}, 
\eqref{estim43} and \eqref{estim preuve lemme carl5} lead to
 \begin{equation*}
 \begin{array}{c}
 I(s,\lambda;\phi)
 \leqslant C\left(s^3\lambda^4
 \displaystyle\iint_{(0,T)\times\omega_2} e^{-2s\alpha}\xi^3|\phi|^2dxdt\right.\\\left.
 ~~~~~~~~~~~~~~~~~~~~~~~~~~~~~~~~~~~~~~~~~~~~~~~~~~~~~~~~~
 +\lambda s^{32/5}(T^2+T^{5/2})\displaystyle\int_{Q_T}e^{-2s\alpha^*}
 (\xi^*)^{32/5}|\psi|^2 dxdt\right),
\end{array}
\end{equation*}
for every $s\geqslant C(T^5+T^{10})$ and $\lambda \geqslant C$. Hence, since 
$$T^2+T^{5/2}\leqslant C s^{2/5}, $$
we have 
 \begin{equation}\label{preuve ine obs 1}
 \begin{array}{c}
 I(s,\lambda;\phi)
 \leqslant C\left(s^3\lambda^4
 \displaystyle\iint_{(0,T)\times\omega_2} e^{-2s\alpha}\xi^3|\phi|^2dxdt
 +\lambda s^{34/5}\displaystyle\int_{Q_T}e^{-2s\alpha^*}
 (\xi^*)^{34/5}|\psi|^2 dxdt\right),
\end{array}
\end{equation}
for every $s\geqslant C(T^5+T^{10})$ and $\lambda \geqslant C$.

\vspace*{0.5cm}

\textbf{Step 2:} We apply Lemma \ref{poincare poids} successively to $\mathcal{N}^*\psi$ with $r=5$, 
then to $\nabla\mathcal{N}^*\psi$ with $r=3$, and we obtain
\begin{equation}\label{preuve ine obs 2}\begin{array}{r}
 s^7\lambda^8\displaystyle\iint_{Q_T} e^{-2s\alpha}\xi^{7}|\mathcal{N}^*\psi|^2dxdt
 \leqslant C\left(
 s^5\lambda^6 \displaystyle\iint_{Q_T} e^{-2s\alpha}\xi^{5}|\nabla \mathcal{N}^*\psi|^2dxdt
 \right.~~~~~~~~~~~~~~~\\\left.
  + s^7\lambda^8\displaystyle\iint_{(0,T)\times\omega_2} e^{-2s\alpha}\xi^{7}|\mathcal{N}^*\psi|^2dxdt\right)
\end{array}\end{equation}
and
\begin{equation}\label{preuve ine obs 3}\begin{array}{r}
 s^5\lambda^6\displaystyle\iint_{Q_T} e^{-2s\alpha}\xi^{5}|\nabla\mathcal{N}^*\psi|^2dxdt
 \leqslant C\left(
 s^3\lambda^4 \displaystyle\iint_{Q_T} e^{-2s\alpha}\xi^{3}|\nabla\nabla \mathcal{N}^*\psi|^2dxdt
 \right.~~~~~~~~~~~~~~~\\\left.
  + s^5\lambda^6\displaystyle\iint_{(0,T)\times\omega_2} e^{-2s\alpha}\xi^{5}|\nabla\mathcal{N}^*\psi|^2dxdt\right),
\end{array}\end{equation}
for every $\lambda\geqslant C$ and $s\geqslant C(T^5+T^{10})$.
 A combination of inequalities \eqref{preuve ine obs 1}-\eqref{preuve ine obs 3} gives
 \begin{equation}\label{est2phi}\begin{array}{l}
\displaystyle\iint_{Q_T}e^{-2s\alpha}\{ 
s^7\lambda^8\xi^{7}|\mathcal{N}^*\psi|^2+s^5\lambda^6\xi^{5}|\nabla\mathcal{N}^*\psi|^2
+s^3\lambda^4\xi^{3}|\nabla\nabla\mathcal{N}^*\psi|^2+s\lambda^2\xi|\nabla\nabla\nabla\mathcal{N}^*\psi|^2   \}dxdt\\
 \leqslant C\left(\lambda s^{34/5}
 \displaystyle\int_{Q_T}e^{-2s\alpha^*}(\xi^*)^{34/5}|\psi|^2 dxdt
 \right.\\\left.
~~~~~~~~~~~~~~~~ +\displaystyle\iint_{(0,T)\times\omega_2}e^{-2s\alpha}\{ 
s^7\lambda^8\xi^{7}|\mathcal{N}^*\psi|^2+s^5\lambda^6\xi^{5}|\nabla\mathcal{N}^*\psi|^2
+s^3\lambda^4\xi^{3}|\nabla\nabla\mathcal{N}^*\psi|^2\}dxdt\right).
 \end{array}\end{equation}

 \textbf{Step 3:} 
Let us consider $\theta_1\in \mathcal{C}^2(\overline{\Omega})$ such that 
\begin{equation*}
 \left\{\begin{array}{ll}
  \Supp(\theta_1)\subseteq \omega_1,&\\
  \theta_1\equiv 1&\mathrm{~in~} \omega_2,\\
  0\leqslant \theta_1\leqslant 1 &\mathrm{~in~} \Omega.
 \end{array}
\right.
\end{equation*}
 Then, after an integration by parts, 
 \begin{equation}\label{step3 1}
 \begin{array}{l}
  s^3\lambda^4\displaystyle\iint_{(0,T)\times\omega_2}e^{-2s\alpha}\xi^{3}|\nabla\nabla\mathcal{N}^*\psi|^2dxdt\\
  \leqslant s^3\lambda^4\displaystyle\iint_{(0,T)\times\omega_1}\theta_1e^{-2s\alpha}\xi^{3}|\nabla\nabla\mathcal{N}^*\psi|^2dxdt\\
   = -s^3\lambda^4\displaystyle\iint_{(0,T)\times\omega_1}
   \sum\limits_{i,j=1}^{N}\{\partial_i(\theta_1e^{-2s\alpha}\xi^{3})\partial_i\partial_j\mathcal{N}^*\psi
   +\theta_1e^{-2s\alpha}\xi^{3}\partial_i^2\partial_j\mathcal{N}^*\psi\}\partial_j(\mathcal{N}^*\psi) dxdt\\
    \leqslant Cs^3\lambda^4\displaystyle\iint_{(0,T)\times\omega_1}
   \{|\nabla(\theta_1e^{-2s\alpha}\xi^{3})||\nabla\nabla\mathcal{N}^*\psi||\nabla\mathcal{N}^*\psi|
   +\theta_1e^{-2s\alpha}\xi^{3}|\nabla\nabla\nabla\mathcal{N}^*\psi||\nabla\mathcal{N}^*\psi|\} dxdt.
   \end{array} \end{equation}
Using the definition of $\xi$ and $\alpha$ given in \eqref{defax}, we deduce that
\begin{equation}\label{estim grad poid}
  |\nabla(\theta_1e^{-2s\alpha}\xi^{3})|\leqslant Cs\lambda e^{-2s\alpha}\xi^{4},
\end{equation}
which, combined with Young's inequality, leads, for every $\varepsilon>0$, to
 \begin{equation}\label{step3 2}\begin{array}{c}
  s^3\lambda^4\displaystyle\iint_{(0,T)\times\omega_2}e^{-2s\alpha}\xi^{3}|\nabla\nabla\mathcal{N}^*\psi|^2dxdt\\
  \leqslant  C\displaystyle\iint_{(0,T)\times\omega_1}e^{-2s\alpha}\{ 
\varepsilon s^3\lambda^4\xi^3|\nabla\nabla\mathcal{N}^*\psi|^2+\varepsilon s\lambda^2\xi|\nabla\nabla\nabla\mathcal{N}^*\psi|^2
+C_{\varepsilon}s^5\lambda^6\xi^{5}|\nabla\mathcal{N}^*\psi|^2\}dxdt,
 \end{array} \end{equation}
 where $C_{\varepsilon}$ depends only on $\varepsilon$. Thus, thanks to inequalities \eqref{est2phi} and \eqref{step3 1},  
   one can absorb (by taking $\varepsilon$ small enough) 
 the local terms  involving $|\nabla\nabla\mathcal{N}^*\psi|^2$ and $|\nabla\nabla\nabla\mathcal{N}^*\psi|^2$ 
 into the right-hand side of inequality 
 \eqref{step3 2} and obtain 
  \begin{equation}\label{step3 3}
  \begin{array}{l}
\displaystyle\iint_{Q_T}e^{-2s\alpha}\{ 
s^7\lambda^8\xi^{7}|\mathcal{N}^*\psi|^2+s^5\lambda^6\xi^{5}|\nabla\mathcal{N}^*\psi|^2
+s^3\lambda^4\xi^{3}|\nabla\nabla\mathcal{N}^*\psi|^2+s\lambda^2\xi|\nabla\nabla\nabla\mathcal{N}^*\psi|^2  \}dxdt\\
 \leqslant C\left(\lambda s^{34/5}
 \displaystyle\int_{Q_T}e^{-2s\alpha^*}(\xi^*)^{34/5}|\psi|^2 dxdt\right.\hspace*{55mm}\\\left.\hspace*{25mm}
 +\displaystyle\iint_{(0,T)\times\omega_1}e^{-2s\alpha}\{ 
s^7\lambda^8\xi^{7}|\mathcal{N}^*\psi|^2+s^5\lambda^6\xi^{5}|\nabla\mathcal{N}^*\psi|^2\}dxdt\right).
 \end{array}\end{equation}
Using exactly the same reasoning we just performed for the term $$s^3\lambda^4\displaystyle\iint_{(0,T)\times\omega_2}e^{-2s\alpha}\xi^{3}|\nabla\nabla\mathcal{N}^*\psi|^2dxdt,$$
 one can also absorb the term $s^5\lambda^6\iint_{(0,T)\times\omega_1}e^{-2s\alpha}\xi^{5}|\nabla\mathcal{N}^*\psi|^2dxdt$ in the right-hand side of \eqref{step3 3} and we obtain
   \begin{equation} \label{ine obs fin 1}\begin{array}{l}
\displaystyle\iint_{Q_T}e^{-2s\alpha}\{ 
s^7\lambda^8\xi^{7}|\mathcal{N}^*\psi|^2+s^5\lambda^6\xi^{5}|\nabla\mathcal{N}^*\psi|^2
+s^3\lambda^4\xi^{3}|\nabla\nabla\mathcal{N}^*\psi|^2+s\lambda^2\xi|\nabla\nabla\nabla\mathcal{N}^*\psi|^2 \}dxdt\\
 ~~~~~~~~~~\leqslant C\left(
\lambda s^{34/5} \displaystyle\int_{Q_T}e^{-2s\alpha^*}(\xi^*)^{34/5}|\psi|^2 dxdt
 +s^7\lambda^8\displaystyle\iint_{(0,T)\times\omega_0}e^{-2s\alpha}
\xi^{7}|\mathcal{N}^*\psi|^2dxdt\right).
 \end{array}\end{equation}
Applying now Lemma \ref{poincare}, and using the definitions of $\alpha^*$ and $\xi^*$ given in \eqref{defaxs}, we obtain the following inequalities:
\begin{equation}\label{ine obs fin 2}
\begin{array}{rcl}
  s^7\lambda^8\displaystyle\iint_{Q_T}(\xi^*)^{7}e^{-2s\alpha^*}|\psi|^2dxdt
  &\leqslant&
  Cs^7\lambda^8\displaystyle\iint_{Q_T}(\xi^*)^{7}e^{-2s\alpha^*} |\mathcal{N}^*\psi|^2dxdt\\
  &\leqslant&
  Cs^7\lambda^8\displaystyle\iint_{Q_T}\xi^{7}e^{-2s\alpha} |\mathcal{N}^*\psi|^2dxdt.
  \end{array}
\end{equation}
The two last inequalities \eqref{ine obs fin 1} and \eqref{ine obs fin 2} give
   \begin{equation*}
   \begin{array}{l}
    s^7\lambda^8\displaystyle\iint_{Q_T}(\xi^*)^{7}e^{-2s\alpha^*}|\psi|^2dxdt\\
~~~~~+\displaystyle\iint_{Q_T}e^{-2s\alpha}\{ 
s^7\lambda^8\xi^{7}|\mathcal{N}^*\psi|^2+s^5\lambda^6\xi^{5}|\nabla\mathcal{N}^*\psi|^2
+s^3\lambda^4\xi^{3}|\nabla\nabla\mathcal{N}^*\psi|^2 +s\lambda^2\xi|\nabla\nabla\nabla\mathcal{N}^*\psi|^2   \}dxdt\\
 ~~~~~~~~~~\leqslant C\left(
\lambda s^{34/5} \displaystyle\int_{Q_T}e^{-2s\alpha^*}(\xi^*)^{34/5}|\psi|^2 dxdt
 +s^7\lambda^8\displaystyle\iint_{(0,T)\times\omega_0}e^{-2s\alpha}
\xi^{7}|\mathcal{N}^*\psi|^2dxdt\right).
 \end{array}\end{equation*}
 Hence, since $34/5<7$, one can absorb the global term of the right-hand side by taking $s$ 
 large enough and obtain inequality \eqref{ine obs 1}.

 \cqfd

Thanks to our Carleman inequality, we can deduce the following observability inequality: 

\begin{prop}\label{prop ine obs 2}
Then for every $\psi^0\in L^2(\Omega)^m$, 
the solution $\psi$ in $W(0,T)^m$ to System \eqref{primal 2x2 constant dirch} satisfies
 \begin{equation}\label{ine obs 3}
 \begin{array}{c}
 \displaystyle\int_{\Omega} 
|\psi(0,x)|^2dx
\leqslant C_{obs}
\displaystyle\iint_{(0,T)\times\omega_0} e^{-2s_0\alpha}\xi^7
|\mathcal{N}^*\psi|^2dxdt,
\end{array} \end{equation}
where $C_{obs}:=Ce^{C(1+T+1/T^5)}$.
\end{prop}

The proof of Proposition \ref{prop ine obs 2} is very classical and is mainly based on dissipation 
estimates and the fact that the weights are bounded from below by some positive constant 
as soon as we are far from $0$ and $T$ (for example on $(T/4,3T/4)$, together with the fact that 
Lemma \ref{poincare} leads to the inequality
\begin{equation*}
\begin{array}{rcl}
  \displaystyle\iint_{(T/4,3T/4)\times\Omega}|\psi|^2dxdt
  &\leqslant&
  C\displaystyle\iint_{(T/4,3T/4)\times\Omega} |\mathcal{N}^*\psi|^2dxdt.
  \end{array}
\end{equation*}
\subsection{Analytic resolution}\label{sec const cont}

 This section is devoted to constructing a solution to the analytic 
control problem \eqref{solution deux controle}, with a control regular enough belonging to the range  
of the differential operator  $\mathcal{N}$. We recall that the definition of $\mathcal{N}$ is given in \eqref{def N}.
Let us consider $\theta\in \mathcal{C}^2(\overline{\Omega})$ such that 
  \begin{equation}\label{deftheta}
 \left\{\begin{array}{ll}
  \Supp(\theta)\subseteq \omega,&\\
  \theta\equiv1&\mathrm{~in~} \omega_0,\\
  0\leqslant \theta\leqslant 1 &\mathrm{~in~} \Omega.
 \end{array}
\right.
\end{equation}

\begin{prop}\label{prop cont reg}Let us assume that Condition \eqref{cond ine} holds. 
Consider the system
 \begin{equation}\label{primal 2x2 contr}
 \left\{\begin{array}{ll}
\partial_tz=\Div(D\nabla z)+G\cdot\nabla  z+Az+\mathcal{N}(\theta v)&\mathrm{in}~ Q_T,\\
       z=0&\mathrm{on}~\Sigma_T,\\
       z(0,\cdot)=y^0&\mathrm{in}~\Omega.
        \end{array}
\right.
\end{equation} 
Then System \eqref{primal 2x2 contr} is null controllable at time $T$, 
i.e. for every $y^0\in L^2(\Omega)^m$ ,
there exists a control $v\in L^2(Q_T)^m$ such that the solution $z$ 
to System \eqref{primal 2x2 contr} satisfies $z(T)\equiv 0$ in $\Omega$. 
Moreover,   for every $K\in(0,1)$, we have 
$e^{Ks_0\alpha^*}v\in W^{2,1}_2(Q_T)^m$ (the definition of $W^{2,1}_2(Q_T)$ is given in \eqref{defW21}) and 
\begin{equation}\label{exp control}
\|e^{Ks_0\alpha^*}v\|_{W^{2,1}_2(Q_T)^m}
\leqslant e^{C(1+T+1/T^5)}\|y^0\|_{L^2(\Omega)^m}.
\end{equation}

\end{prop}

\noindent\textbf{Proof of Proposition \ref{prop cont reg}.}

We will use the usual duality method developed by  Fursikov and  Imanuvilov 
in \cite{fursikov1996controllability} in the spirit of what was done in \cite{MR1751309} to obtain more regular controls. 
Let $y^0\in L^2(\Omega)^m$ 
and $\rho$ be the weight defined by
\begin{equation*}
 \rho:=\xi^{7}e^{-2s_0\alpha}.
\end{equation*}
Let $k\in\mathbb{N}^*$ and let us consider the following optimal control problem
\begin{equation}\label{probleme minimisation}
 \left\{\begin{array}{l}
         \mathrm{minimize}~ J_k(v):=\dfrac12\displaystyle\int_{Q_T}\rho^{-1} |v|^2dxdt
 +\dfrac{k}{2}\displaystyle\int_{\Omega}| z(T)|^2dx,\\
 v\in L^2(Q_T,\rho^{-1/2})^m,
        \end{array}
\right.
\end{equation}
where  $z$ is the solution in $W(0,T)^m$ to System \eqref{primal 2x2 contr}.

The functional $J_k:L^2(Q_T,\rho^{-1/2})^m\rightarrow\mathbb{R}^+$ is differentiable, coercive and strictly 
convex on the space $L^2(Q_T,\rho^{-1/2})^m$. 
Therefore 
there exists a unique solution 
to the control optimal problem \eqref{probleme minimisation} (see \cite[p. 128]{lionscontrole}) 
and the optimal control $v_k$ 
is characterized thanks to the solution $z_k$ of the primal system 
\begin{equation}\label{système ramené a tout l espace}
 \left\{\begin{array}{ll}
       \partial_tz_k=\Div(D\nabla z)_k+G\cdot \nabla z_k+Az_k
       +\mathcal{N}({\theta} v_k)&\mathrm{in}~ Q_T,\\
       z_k=0&\mathrm{on}~\Sigma_T,\\
       z_k(0,\cdot)=y^0&\mathrm{in}~\Omega,
        \end{array}
\right.
\end{equation}
the solution $\varphi_k$ to the dual system
\begin{equation}\label{system dual fursikov}
\left\{\begin{array}{ll}
         -\partial_t\varphi_k=\Div(D\nabla\varphi)_k-G^*\cdot\nabla \varphi_k+A^*\varphi_k&\mathrm{in}~ Q_T,\\
                \varphi_k=0&\mathrm{on}~\Sigma_T,\\
       \varphi_k(T,\cdot)=k  z_k(T,\cdot)&\mathrm{in}~\Omega
        \end{array}
\right.
\end{equation}
and the relation
\begin{equation}\label{carac controle}
\left\{\begin{array}{l}
        v_k=-\rho \theta\mathcal{N}^*\varphi_k\mathrm{~in~}Q_T,\\
       v_k\in L^2(Q_T,\rho^{-1/2})^m.
        \end{array}
\right.
\end{equation}

The rest of the proof is divided into two steps. In the first step, we will prove that the sequence $(v_k)_{k\in \mathbb{N}^*}$ 
converges to a control $v\in L^2(Q_T,\rho^{-1/2})^m$ with an associated solution $z$ to System \eqref{primal 2x2 contr} satisfying  
$z(T)\equiv 0$ in $\Omega$.  Then, in the second step, we will establish \eqref{exp control}.

\vspace*{0.5cm}

\textbf{Step 1:} \\
Firstly, the characterization  (\ref{système ramené a tout l espace}), (\ref{system dual fursikov}) 
and \eqref{carac controle} of the minimizer $v_k$ of $J_k$  in $L^2(Q_T,\rho^{-1/2})^m$ 
leads to the following computations
\begin{equation}\label{estim J_k1}
\begin{array}{rcl}
 J_k(v_k)&=&-\dfrac12\displaystyle\int_{0}^T\langle \theta \mathcal{N}^*\varphi_k,v_k\rangle_{L^2(\Omega)^m}dt
 +\dfrac12\langle z_k(T),\varphi_{k}(T)\rangle_{L^2(\Omega)^m}\vspace*{3mm}\\
 &=&-\dfrac12\displaystyle\int_{0}^T\langle \varphi_k,\mathcal{N}(\theta v_k)\rangle_{L^2(\Omega)^m}dt
  +\dfrac12\displaystyle\int_{0}^T\{\langle z_k,\partial_t\varphi_k\rangle_{L^2(\Omega)^m}
  +\langle\partial_tz_k,\varphi_k\rangle_{L^2(\Omega)^m}\}dt\\
  &&~~~~~~~~~~~~~~~~~~~~~~~~~~~~~~~~~~~~~~~~~~~~~~~~~~~~~~~~~~~~~
  +\dfrac12\langle y^0,\varphi_{k}(0,\cdot)\rangle_{L^2(\Omega)^m}\\
 &=& \dfrac12\langle y^0,\varphi_{k}(0,\cdot)\rangle_{L^2(\Omega)^m}.
 \end{array}
 \end{equation}
 Moreover, using the definition of $J_k$, the definition of $\theta$, our observability inequality \eqref{ine obs 3}, the expression \eqref{estim J_k1} and the Cauchy-Schwarz inequality, we infer 
  \begin{equation*}
  \begin{array}{cl}
 \|\varphi_k(0,\cdot)\|_{L^2(\Omega)^m}^2&\leqslant 
 C_{obs} \displaystyle\iint_{Q_T}\rho\theta^2|\mathcal{N}^*\varphi_k|^2dxdt
 =C_{obs} \displaystyle\iint_{Q_T}\rho^{-1}|v_k|^2dxdt\\
 &\leqslant 2C_{obs}J_k(v_k)\leqslant 2C_{obs} \|\varphi_k(0,\cdot)\|_{L^2(\Omega)^m}\|y^0\|_{L^2(\Omega)^m},
\end{array}
 \end{equation*}
from which we deduce 
  \begin{equation}\label{estim J_k2}
 \|\varphi_k(0,\cdot)\|_{L^2(\Omega)^m}\leqslant 2C_{obs}\|y^0\|_{L^2(\Omega)^m}.
 \end{equation}
 Then, using  \eqref{estim J_k1} and \eqref{estim J_k2}, we deduce 

\begin{equation}\label{J_k(V) borne}
 J_k(v_k)\leqslant{C_{obs}}\|y^0\|_{L^2(\Omega)^m}^2.
\end{equation}
Furthermore, we have (see \cite{lionscontrole})
\begin{equation}\label{norme z bornee}
\begin{array}{rl}
 \|z_k\|_{W(0,T)^m}
 &\leqslant C\left(\|
 \mathcal{N}(\theta v_k)\|_{L^2(0,T;H^{-1}(\Omega))^m}+\|y^0\|_{L^2(\Omega)^m}\right),\\
  &\leqslant C\left(\|
 v_k\|_{L^2(Q_T)}+\|y^0\|_{L^2(\Omega)^m}\right),\\
  &\leqslant C(1+C_{obs})
  \|y^0\|_{L^2(\Omega)^m},
\end{array}
\end{equation}
where $C$ does not depend on $y^0$ and $k$. 
Then, using  inequalities \eqref{J_k(V) borne} and \eqref{norme z bornee}, we deduce 
that there exist subsequences, which are still denoted $v_k$, $z_k$, 
such that the following weak convergences hold:
\begin{equation*}\left\{
 \begin{array}{ll}
  v_k\rightharpoonup v~&\mathrm{in}~L^2(Q_T,\rho^{-1/2})^m,\\
  z_k\rightharpoonup z~&\mathrm{in}~W(0,T)^m,\\
    z_k(T)\rightharpoonup   0~&\mathrm{in}~L^2(\Omega)^m.\\
 \end{array}
\right.
\end{equation*}
Passing to the limit in $k$, $z$ is solution to System (\ref{primal 2x2 contr}). 
Moreover, using the expression of $J_k$ given in \eqref{probleme minimisation} and inequality \eqref{J_k(V) borne}, we deduce by letting $k$ going to $\infty$ that $ z(T)\equiv0$ in $\Omega$. 
Thus the solution $z$ to System \eqref{primal 2x2 contr} with control $v\in L^2(Q_T,\rho^{-1/2})^m$ 
satisfies $z(T)\equiv0$ in $\Omega$ 
and using (\ref{J_k(V) borne}), we obtain the inequality 
\begin{equation}\label{estim contr}
 \|v\|_{L^2(Q_T,\rho^{-1/2})^m}^2\leqslant C_{obs} \|y^0\|_{L^2(\Omega)^m}^2.
\end{equation}

\textbf{Step 2:}  
One remarks that for every $K\in(0,1)$, there exists a constant $C>0$ such that $$e^{2Ks_0\alpha^*}\leqslant C\xi^{-7}e^{2s_0\alpha}.$$ 
This inequality and  estimate (\ref{J_k(V) borne})   imply that
\begin{equation}\label{estim H2 u 0}
  \|e^{Ks_0\alpha^*} v_k\|_{L^2(Q_T)^m}^2\leqslant C\displaystyle\int_{Q_T}\xi^{-7}e^{2s_0\alpha} |v_k|^2dxdt
  \leqslant J_k(v_k)\leqslant e^{C(1+T+1/T^5)}\|y^0\|_{L^2(\Omega)^m}^2.
\end{equation}
We recall that $v_k$  is defined in \eqref{carac controle}, moreover, thanks to the definitions of $\xi$ and $\alpha$ given in \eqref{defax}, one has, for every $\eta\geqslant0$ ,
\begin{equation*}
 \begin{array}{l}
    |\nabla (\xi^{\eta}e^{-2s_0\alpha})|\leqslant C \xi^{\eta+1}e^{-2s_0\alpha},\\
  |\Delta (\xi^{\eta}e^{-2s_0\alpha})|\leqslant C \xi^{\eta+2}e^{-2s_0\alpha},  \\
  |\partial_t (\xi^{\eta}e^{-2s_0\alpha})|\leqslant C T\xi^{\eta+6/5}e^{-2s_0\alpha}.
 \end{array}
\end{equation*}
These above inequalities, for $\eta:=7$, lead to the fact that

\begin{equation}\label{estim H2 u 1}
 \|e^{Ks_0\alpha^*}\nabla v_k\|_{L^2(Q_T)^m}^2\leqslant 
 C\displaystyle\iint_{Q_T}e^{-4s_0\alpha+2Ks_0\alpha^*}\{\xi^{14}|\nabla\mathcal{N}^*\varphi_k|^2+\xi^{16}|\mathcal{N}^*\varphi_k|^2\}dxdt,
\end{equation}
\begin{equation}\label{estim H2 u 2}
 \|e^{Ks_0\alpha^*}\Delta v_k\|_{L^2(Q_T)^m}^2\leqslant 
 C\displaystyle\iint_{Q_T}e^{-4s_0\alpha+2Ks_0\alpha^*}\{\xi^{14}|\nabla\nabla\mathcal{N}^*\varphi_k|^2+\xi^{16}|\nabla\mathcal{N}^*\varphi_k|^2+\xi^{18}|\mathcal{N}^*\varphi_k|^2\}dxdt
\end{equation}
and
\begin{equation}\label{estim H2 u 3}
\begin{array}{c}
 \|\partial_t (e^{Ks_0\alpha^*}v_k)\|_{L^2(Q_T)^m}^2\leqslant CT\displaystyle\iint_{Q_T}e^{-4s_0\alpha+2Ks_0\alpha^*}\{\xi^{14}|\mathcal{N}^*\partial_t\varphi_k|^2+\xi^{82/5}|\mathcal{N}^*\varphi_k|^2\}dxdt\\
 \leqslant CT\displaystyle\iint_{Q_T}e^{-4s_0\alpha+2Ks_0\alpha^*}\{\xi^{14}(|\mathcal{N}^*\varphi_k|^2
 +|\Delta\mathcal{N}^*\varphi_k|^2+|\nabla\mathcal{N}^*\varphi_k|^2)+\xi^{82/5}|\mathcal{N}^*\varphi_k|^2\}dxdt.
\end{array}\end{equation}
For every $\eta,\nu>0$ there exists a constant $C_{\eta,\nu}$ such that
\begin{equation*}
 \begin{array}{l}
  |\xi^{\eta}e^{-4s\alpha+2Ks_0\alpha^*}|\leqslant C_{\eta,\nu} \xi^{\nu}e^{-2s\alpha}.
 \end{array}
\end{equation*}
Combining the last  inequality with \eqref{estim H2 u 0}-\eqref{estim H2 u 3}, 
we deduce that 
\begin{equation*} \begin{array}{c}
 \| e^{Ks_0\alpha^*}v_k\|_{W^{2,1}_2(Q_T)^m}^2\\
 \leqslant e^{C(1+T+1/T^5)}\displaystyle\iint_{Q_T}e^{-2s_0\alpha}\{ 
\xi^{7}|\mathcal{N}^*\varphi_k|^2+\xi^{5}|\nabla\mathcal{N}^*\varphi_k|^2
+\xi^{3}|\nabla\nabla\mathcal{N}^*\varphi_k|^2dxdt.
 \end{array}\end{equation*}
Using \eqref{ine obs 1}, we obtain that  $e^{Ks_0\alpha^*}\nabla v_k$, 
$e^{Ks_0\alpha^*}\Delta v_k$, $\partial_t (e^{Ks_0\alpha^*}v_k)\in L^2(Q_T)^m$, and that
\begin{equation*}
 \| e^{Ks_0\alpha^*}v_k\|_{W^{2,1}_2(Q_T)^m}^2
 \leqslant e^{C(1+T+1/T^5)}\displaystyle\iint_{Q_T} e^{-2s\alpha}\xi^7
|\theta\mathcal{N}^*\varphi_k|^2dxdt=e^{C(1+T+1/T^5)}\|v_k\|_{L^2(Q_T)^m}^2.
\end{equation*}
The estimate \eqref{J_k(V) borne} of $v_k$ gives
\begin{equation*}
 \|e^{Ks_0\alpha^*}v_k\|_{W^{2,1}_2(Q_T)^m}\leqslant e^{C(1+T+1/T^5)}\|y^0\|_{L^2(\Omega)^m}.
\end{equation*}
We conclude by letting $k\rightarrow+\infty$ (after extracting an adequate subsequence) in the inequalities above.

\cqfd

\subsection{End of the proof of Theorem \ref{th 1}}\label{sec proof th1}

Let us assume that Condition \ref{cond th 1} holds. We will prove first  
the null controllability at time $T$ of System \eqref{system primmal}. 
 Let $y^0\in L^2(\Omega)^m$. 
We follow the method explained in Section \ref{section strategy}. Let us remind that $\theta$ is defined in \eqref{deftheta}.
Using  Proposition \ref{prop cont reg},
 the following system:
 \begin{equation}\label{primal 2x2 contr th1}
 \left\{\begin{array}{ll}
\partial_tz=\Div(D\nabla z)+G\cdot\nabla  z+Az+\mathcal{N}(\theta v)&\mathrm{in}~ Q_T,\\
       z=0&\mathrm{on}~\Sigma_T,\\
       z(0,\cdot)=y^0&\mathrm{in}~\Omega
        \end{array}
\right.
\end{equation} 
 is null controllable at time $T$,  
thus there exists a control $v\in L^2(Q_T)^m$ such that the solution $z$ in $W(0,T)^m$ 
to System \eqref{primal 2x2 contr th1} satisfies 
\begin{equation*}
z(T)\equiv 0 \mathrm{~~in~~} \Omega.  
\end{equation*}
Moreover  
\begin{equation}\label{exp control th1}
 e^{Ks_0\alpha^*}v\in W^{2,1}_2(Q_T)^m.
\end{equation}
Taking into account Proposition \ref{alg1} and the definition of $\mathcal M^*$ given in \eqref{def M*}, one has \eqref{LM=N} with the operator
\begin{equation*}
\begin{array}{cccc}
 \mathcal M:&W^{2,1}_2(Q_T)^m&\rightarrow &L^2(Q_T)^{m}\times L^2(Q_T)^{m-1}\\
 &f&\mapsto &\mathcal{M}f,
\end{array}\end{equation*}
defined by 
\begin{equation*}
 \mathcal{M}f=\left(\begin{array}{c}
       0\\\vdots\\0\\\phantom{~~~~(i_0^{th}~\mathrm{line})} f_m~~~~(i_0^{th}~\mathrm{line})\\0\\\vdots\\0\\
       -(g_{mi_0}\cdot\nabla+a_{mi_0})f_1+(g_{1i_0}\cdot\nabla+a_{1i_0})f_m\\\vdots\\
       -(g_{mi_0}\cdot\nabla+a_{mi_0})f_{i_0-1}+(g_{(i_0-1)i_0}\cdot\nabla+a_{(i_0-1)i_0})f_{m}\\
       (-\partial_t+\Div(d_{i_0}\nabla f_m) -(g_{mi_0}\cdot\nabla+a_{mi_0})f_{i_0}+(g_{i_0i_0}\cdot\nabla+a_{i_0i_0})f_{m}\\
       -(g_{mi_0}\cdot\nabla+a_{mi_0})f_{i_0+1}+(g_{(i_0+1)i_0}\cdot\nabla+a_{(i_0+1)i_0})f_{m}\\\vdots\\
       -(g_{mi_0}\cdot\nabla+a_{mi_0})f_{m-1}+(g_{(m-1)i_0}\cdot\nabla+a_{(m-1)i_0})f_{m}
                \end{array}\right).
\end{equation*}
Let $(\widehat{z},\widehat{v})$ be defined by
\begin{equation*}
\left(\begin{array}{c}\widehat{z}\\\widehat{v}\end{array}\right)
 :=\mathcal{M}\left(\theta v\right).
\end{equation*}
Using \eqref{exp control th1} and the fact that $\mathcal M$ is a differential operator of order $1$ in time and $2$ in space (see Proposition \ref{alg1}) 
with bounded coefficients, we obtain that $(\widehat{z},\widehat{v})\in L^2(Q_T)^m\times L^2(Q_T)^{m-1}$. 
Moreover, using \eqref{exp control th1}, we have $\widehat{z}(0,\cdot)=\widehat{z}(T,\cdot)=0$ in $\Omega$ and 
we remark that $(\widehat{z},\widehat{v})$ 
is a solution 
to the control problem
\begin{equation}\label{probleme ramene a tout l espace th1}
 \left\{\begin{array}{ll}
       \partial_t\widehat{z}=\Div(D\nabla  \widehat{z})+G\cdot\nabla \widehat{z}+A\widehat{z}+B\widehat{v}+\mathcal{N}(\theta v)&\mathrm{in}~ Q_T,\\
       \widehat{z}=0&\mathrm{on}~\Sigma_T,\\
       \widehat{z}(0,\cdot)=\widehat{z}(T,\cdot)=0&\mathrm{in}~\Omega, 
        \end{array}
\right.
\end{equation}
in particular $\mathcal{L}\circ\mathcal{M}=\mathcal{N}$. 
Finally, $(\widehat{z},\widehat{v})\in W(0,T)^m\times L^2(Q_T)^{m-1}$ 
thanks to the usual parabolic regularity. 
Thus the pair $(y,u):=(z-\widehat{z},-\widehat{v})$ is a solution to System \eqref{system primmal} 
in $W(0,T)^m\times L^2(Q_T)^{m-1}$ and satisfies 
\begin{equation*}
y(T,\cdot)\equiv 0 \mathrm{~~in~~} \Omega. 
\end{equation*}

\section{Proof of Theorem \ref{th 2}}
Let us remind that in this case, we have only $2$ equations.

\subsection{Algebraic resolution}

We will assume in this section that all differential operators are defined in $\mathcal{C}^{\infty}(Q_T)^2$. 
We recall that  $\mathcal{N}$ is simply the identity operator and $\mathcal{L}$ is given in \eqref{def L Strategy}.  
We want to find a differential operator $\mathcal{M}$ satisfying
\begin{equation}\label{LM=I}
 \mathcal{L}\circ\mathcal{M}=Id.
\end{equation}

Let us emphasize that when coefficients are depending on time and space, we prove equality \eqref{LMB} with 
the identity operator in the right-hand side (and not $\mathcal N$ as defined in \eqref{def N}), 
that is equality \eqref{LM=I}, because Proposition  
\ref{prop cont reg} holds only for constant coefficients and does not seem to be adapted to the case of non-constant coefficients.
We have the following proposition:
 \begin{prop}\label{alg2} 
 One has
 \begin{enumerate}
  \item[(i)] Under Conditions \eqref{cond1th2} and \eqref{reg1th2}, there exists a differential operator 
        $\mathcal{M}$ of order at most $1$ in time and  $2$ in space, 
        with bounded coefficients on $(a,b)\times{\mathcal O}$, such that  equality \eqref{LM=I} holds.
  \item[(ii)]If $N=1$ and under Condition \eqref{cond2th2} and \eqref{reg2th2}, there exists a differential operator 
        $\mathcal{M}$ of order at most $2$ in time,  $4$ in space, and $1-2$ respectively in crossed space-time, 
        with bounded coefficients on $(a,b)\times{\mathcal O}$, such that  equality \eqref{LM=I} holds.
 \end{enumerate}

\end{prop}
\noindent\textbf{Proof of Proposition \ref{alg2}.}

Equality \eqref{LM=I} is equivalent to
\begin{equation}\label{M*L*=I}
 \mathcal{M}^*\circ\mathcal{L^*}=Id.
\end{equation}
Taking into account the definition of $\mathcal L$ given in \eqref{def L Strategy}, 
the adjoint $\mathcal L^*$ of the operator $\mathcal{L}$ is given by
\begin{equation}\label{def L* 3x3}
 \mathcal{L}^*\varphi:= \left(\begin{array}{c}
\mathcal{L}^*_1\varphi\\
\mathcal{L}^*_2\varphi\\
\mathcal{L}^*_3\varphi
        \end{array}\right)=
 \left(\begin{array}{c}
-\partial_t\varphi_1-\Div(d_1\nabla \varphi_1)
+\Div(g_{11} \varphi_1)+\Div (g_{21}\varphi_2)-a_{11}\varphi_1-a_{21}\varphi_2\\
-\partial_t\varphi_2-\Div(d_2\nabla \varphi_2)
+\Div (g_{12} \varphi_1)+\Div( g_{22}\varphi_2)-a_{12}\varphi_1-a_{22}\varphi_2\\
\varphi_1
        \end{array}\right).
\end{equation}
We remark first that 
\begin{equation}\label{base2}\mathcal{L}^*_1\varphi
+\{\partial_t+\Div(d_1\nabla\cdot)-\Div(g_{11}\cdot)+a_{11}\}\circ\mathcal{L}^*_3\varphi
=\Div(g_{21}\varphi_2)-a_{21}\varphi_2.\end{equation}

   \begin{enumerate}
    \item[(i)]
Let us consider some open subset $\widetilde {\mathcal O}$  included in $\mathcal O$ on which we have $|a_{21}|\geqslant C>0$ (such an open subset exists thanks to \eqref{cond1th2}).
    Since $g_{21}=0$ in $(a,b)\times{\mathcal O}$, one can just consider  $\mathcal{M}^*$ defined  for every $\psi:=(\psi_1,\psi_2,\psi_3)\in\mathcal{C}^{\infty}(Q_T)^3$, 
locally on $(a,b)\times{\widetilde {\mathcal O}}$, by
\begin{equation}\label{M*21}
 \mathcal{M}^* \psi
        :=\left(\begin{array}{c}
\psi_3\\
-\frac{\psi_1
+\{\partial_t+\Div(d_1\nabla\cdot)-\Div(g_{11}\cdot)+a_{11}\}\psi_3}{a_{21}}
        \end{array}\right),
        \end{equation}
        so that equality \eqref{M*L*=I} is satisfied. Moreover, the coefficients of $\mathcal M^*$ (and hence of $\mathcal M$) are bounded.

        \item[(ii)]If $N=1$ and under Condition \eqref{cond2th2}, 
     one can proceed as follows: we consider the operator
      \begin{equation}\label{defQ}\mathcal Q(\varphi)=\left(\begin{array}{c}
        \mathcal{L}^*_3\varphi\\
      \mathcal{L}^*_1\varphi
+\{\partial_t+\partial_x(d_1\partial_x\cdot)-\partial_x(g_{11}\cdot)+a_{11}\}\circ\mathcal{L}^*_3\varphi
\\\partial_x(
\mathcal{L}^*_1\varphi
+\{\partial_t+\partial_x(d_1\partial_x\cdot)-\partial_x(g_{11}\cdot)+a_{11}\}\circ\mathcal{L}^*_3\varphi)
\\\partial_t(\mathcal{L}^*_1\varphi
+\{\partial_t+\partial_x(d_1\partial_x\cdot)-\partial_x(g_{11}\cdot)+a_{11}\}\circ\mathcal{L}^*_3\varphi)\\
\partial_{xx}(\mathcal{L}^*_1\varphi
+\{\partial_t+\partial_x(d_1\partial_x\cdot)-\partial_x(g_{11}\cdot)+a_{11}\}\circ\mathcal{L}^*_3\varphi)\\
\mathcal{L}^*_2\varphi+\{a_{12}-\partial_x(g_{12}\cdot)\}\circ\mathcal{L}^*_3\varphi,\\
\partial_x(\mathcal{L}^*_2\varphi+\{a_{12}-\partial_x(g_{12}\cdot)\}\circ\mathcal{L}^*_3\varphi)
            \end{array}\right ),
      \end{equation}
      i.e.
            \begingroup\scriptsize\begin{equation}\label{defQ2}
            \mathcal Q(\varphi)=\left(\begin{array}{c}
 \varphi_1\\
(-a_{21}+\partial_x g_{21})\varphi_2+g_{21}\partial_x \varphi_2\\
(-\partial_x a_{21}+\partial_{xx}g_{21})\varphi_2+(-a_{21}+2\partial_x g_{21})\partial_x\varphi_2
+g_{21}\partial_{xx} \varphi_2\\
(-\partial_t a_{21}+\partial_{tx} g_{21})\varphi_2+\partial_t g_{21}\partial_x\varphi_2
+(-a_{21}+\partial_{x} g_{21})\partial_t\varphi_2
+ g_{21}\partial_{xt} \varphi_2\\
(-\partial_{xx} a_{21}+\partial_{xxx}g_{21}) \varphi_2+(-2\partial_{x} a_{21}
+3\partial_{xx} g_{21})\partial_x\varphi_2+(-a_{21}+3\partial_{x} g_{21})\partial_{xx}\varphi_2
+g_{21}\partial_{xxx} \varphi_2\\
(-a_{22}+\partial_x g_{22})\varphi_2+(-\partial_x d_2+g_{22})\partial_x\varphi_2
-\partial_t\varphi_2-d_2\partial_{xx}\varphi_2\\
(-\partial_x a_{22}+\partial_{xx} g_{22})\varphi_2
+(-\partial_{xx} d_2-a_{22}+2\partial_x g_{22})\partial_x\varphi_2
+(-2\partial_x d_2+g_{22})\partial_{xx}\varphi_2-\partial_{tx}\varphi_2-d_2\partial_{xxx}\varphi_2
            \end{array}\right ).
      \end{equation}\endgroup
      
It is easy to see that there are only $7$ different derivatives of $\varphi$ appearing in \eqref{defQ}, 
which are the following ones: $$\varphi_1,\varphi_2,\partial_{x} \varphi_2,\partial_{t} \varphi_2,
\partial_{xx} \varphi_2,\partial_{xt} \varphi_2,\partial_{xxx} \varphi_2.$$
 Hence, we can see the operator $\mathcal Q$  as a matrix $M$ acting on these derivatives 
 (see \cite[Sec. 3.2]{coronlissy2014}), so that $M$ is a square matrix of size $7\times 7$. 
 More precisely, one has 
 \begin{equation}\label{decomp Q avec M}
 \mathcal Q(\varphi)=M(\varphi_1,\varphi_2,\partial_{x} \varphi_2,\partial_{t} \varphi_2
  ,\partial_{xx} \varphi_2,\partial_{xt} \varphi_2,\partial_{xxx} \varphi_2),
  \end{equation}
 with
  \begingroup\scriptsize
  \begin{equation}\label{def M alg}M:=\left(\begin{array}{cccccccc} 1&0&0&0&0&0&0\\
0&-a_{21}+\partial_x g_{21}&g_{21}&0&0&0&0\\
0&-\partial_x a_{21}+\partial_{xx} g_{21}&-a_{21}+2\partial_x g_{21}&0&g_{21}&0&0\\
0&-\partial_t a_{21}+\partial_{tx} g_{21}&\partial_{t} g_{21}&-a_{21}+\partial_x g_{21}&0&g_{21}&0\\
0&-\partial_{xx} a_{21}+\partial_{xxx} g_{21}&-2\partial_x a_{21}+3\partial_{xx} g_{21}&0&
- a_{21}+3\partial_{x} g_{21}&0&g_{21}\\
0&-a_{22}+\partial_xg_{22}&-\partial_xd_2+g_{22}&-1&-d_2&0&0\\
0&-\partial_x a_{22}+\partial_{xx}g_{22}&-\partial_{xx}d_2-a_{22}+2\partial_xg_{22}&0
&-2\partial_xd_2+g_{22}&-1&-d_2\end{array}\right ).
  \end{equation}\endgroup
Matrix $M$ is invertible since Condition \eqref{cond2th2} is verified.
Let us call $P$ the projection on the two first components 
$$P(x_1,x_2,x_3,x_4,x_5,x_6,x_7)=(x_1,x_2).$$
Then, by definition of the inverse, we have
$$PM^{-1}M (\varphi_1,\ldots,\partial_{xxx}\varphi_2)=\varphi.$$
The expression of $\mathcal{Q}$ given in \eqref{decomp Q avec M} leads to 
$$PM^{-1}\mathcal Q(\varphi)
=\varphi.$$
If we denote by $\mathcal{R}_1:=\partial_t+\partial_x(d_1\partial_x\cdot)-\partial_x(g_{11}\cdot)+a_{11}$ 
and $\mathcal{R}_2:=a_{12}-\partial_x(g_{12}\cdot)$, 
using \eqref{defQ}, we remark that the previous equality can be rewritten as 
$$
PM^{-1}\mathcal{S}\circ\mathcal{L}^*\varphi
=\varphi,$$
where 
\begin{equation}\label{def S alg}
\mathcal{S}:=\left(\begin{array}{ccc}0&0&1\\1&0&\mathcal{R}_1\\\partial_x&0&\partial_x\circ\mathcal{R}_1\\
\partial_t&0&\partial_t\circ\mathcal{R}_1\\\partial_{xx}&0&\partial_{xx}\circ\mathcal{R}_1\\0&1&\mathcal{R}_2\\
0&\partial_x&\partial_x \circ\mathcal{R}_2\end{array}\right).
\end{equation}
Hence equality \eqref{M*L*=I} is satisfied for 
\begin{equation}\label{M*22}
 \mathcal{M}^*:=PM^{-1}\mathcal{S}.
\end{equation}
Moreover, thanks to Conditions \eqref{cond2th2},  the coefficients of $\mathcal M^*$ (and hence of $\mathcal M$) are bounded.
   \end{enumerate}

\cqfd

\begin{rem}
The most important point in the construction of $\mathcal Q$ is to differentiate enough time the equations as in \eqref{defQ} in order to obtain the same number equations than ``unknowns'', the unknowns being $\varphi$ and its derivatives seen as independent algebraic variables (see notably \cite[Section 3.2.2]{coronlissy2014} for further explanations). One can check that this cannot be performed by differentiating the equations less times.

\end{rem}

\subsection{Analytic resolution}
Let $\omega_1$ be a nonempty open subsets included in $\omega$ satisfying
\begin{equation*}
 \overline{\omega}_1\subset\omega_0.
\end{equation*}
Let us consider $\theta\in \mathcal{C}^2(\overline{\Omega})$ such that 
  \begin{equation}\label{theta2}
 \left\{\begin{array}{ll}
  \Supp(\theta)\subseteq \omega_0,&\\
  \theta\equiv1&\mathrm{~in~} \omega_1,\\
  0\leqslant \theta\leqslant 1 &\mathrm{~in~} \Omega.
 \end{array}
\right.
\end{equation}
We are going to explain what are the main differences with Subsection \ref{sec ine obs} in order to obtain a Carleman inequality.
First of all, we need to find an equivalent to Lemma \ref{th carl neum}, which is the following:

\begin{Lemme}\label{lemme 3.1}
 Let us assume that $u^0\in L^2(\Omega)$, $f_1\in L^2(Q_T)$ and $d\in W^{1}_{\infty}(Q_T)$ 
 such that $d>C>0$ in $Q_T$.  Then there exists a constant $C:=C(\Omega,\omega_2)>0$ such that the solution to the system 
 \begin{equation*}\left\{\begin{array}{ll}
  -\partial_tu-\Div(d\nabla u) =f_1&\mathrm{in}~Q_T,\\
  \frac{\partial u}{\partial n} =f_2&\mathrm{on}~\Sigma_T,\\
  u(T,\cdot)=u^0&\mathrm{in~}\Omega,\end{array}\right.
 \end{equation*}
satisfies
\begin{equation*}\begin{array}{r}
 I(s,\lambda;u)\leqslant C\left(
 s^3\lambda^4\displaystyle\iint_{(0,T)\times\omega_2} e^{-2s\alpha}\xi^3|u|^2dxdt
 +\displaystyle\iint_{Q_T} e^{-2s\alpha}|f_1|^2dxdt
 ~~~~~~~~~~~~~~~~~~~~~\right.\\\left.
 +s\lambda\displaystyle\iint_{\Sigma_T}e^{-2s\alpha^*}\xi^*|f_2|^2 d\sigma dt\right),
\end{array}\end{equation*}
for all $\lambda\geqslant C$ and $s\geqslant C(T^5+T^{10})$.
\end{Lemme}
The proof of Lemma \ref{lemme 3.1} can be easily obtained by using the method of \cite{GuerreroFourier} together with the Carleman estimate proved in \cite{yamamoto2003carleman}. 
Let us consider the backward system
\begin{equation}\label{carlpb2}
 \left\{\begin{array}{ll}
-\partial_t\psi_1=\Div(d_1\nabla \psi_1)
-\Div (g_{11}\psi_1)-\Div(g_{21}\psi_2)+a_{11}\psi_1+a_{21}\psi_2&\mathrm{in}~ Q_T,\\
-\partial_t\psi_2=\Div(d_2\nabla \psi_2)
-\Div (g_{12}\psi_1)-\Div (g_{22}\psi_2)+a_{12}\psi_1+a_{22}\psi_2&\mathrm{in}~ Q_T,\\
       \psi_1=\psi_2=0&\mathrm{on}~\Sigma_T,\\
       \psi_1(T,\cdot)=\psi_1^0,~\psi_2(T,\cdot)=\psi_2^0&\mathrm{in}~\Omega,
        \end{array}
\right.
\end{equation}
where $\psi^0:=(\psi_1^0,\psi_2^0)\in L^2(\Omega)^2$. 
From the last Lemma, one can deduce:

\begin{prop}\label{prop ine obs esp temps}
One has 
\begin{enumerate}
 \item[(i)]Under Condition \eqref{cond1th2} and \eqref{reg1th2}, 
 there exists a   constant $C:=C(\omega_1,\Omega)>0$ such that for every 
  $\psi^0:=(\psi_1^0,\psi_2^0)\in L^2(\Omega)^2$, 
 the corresponding solution 
 $\psi:=(\psi_1,\psi_2)$  of the backward problem \eqref{carlpb2} 
satisfies
\begin{equation*}\label{derniercarl1}
 \begin{array}{c}
\displaystyle\iint_{Q_T} e^{-2s\alpha}\{s^5\lambda^6\xi^{5}|\psi|^2
+s^3\lambda^4\xi^{3}|\nabla\psi|^2+s\lambda^2\xi|\nabla\nabla\psi|^2\}dxdt\\
\leqslant C s^5\lambda^{6}
\displaystyle\iint_{(0,T)\times\omega_1} e^{-2s\alpha}\xi^5
|\psi|^2dxdt.
\end{array} 
\end{equation*}
 \item[(ii)] If $N=1$ and under Condition  \eqref{cond2th2} and \eqref{reg2th2}, we obtain the same conclusion as in item (i) by replacing estimate \eqref{derniercarl1} by
\begin{equation}
 \begin{array}{c}
\displaystyle\iint_{Q_T} e^{-2s\alpha}\{s^9\lambda^{10}\xi^{9}|\psi|^2+
s^7\lambda^8\xi^{7}|\nabla\psi|^2+s^5\lambda^6\xi^{5}|\nabla\nabla\psi|^2
\\+s^3\lambda^4\xi^{3}|\nabla\nabla\nabla\psi|^2+s\lambda^2\xi|\nabla\nabla\nabla\nabla\psi|^2\\
\leqslant C s^9\lambda^{10}
\displaystyle\iint_{(0,T)\times\omega_1} e^{-2s\alpha}\xi^9
|\psi|^2dxdt.
\end{array} 
\end{equation}
\end{enumerate} 
\end{prop}

The proof is very similar to the proof of Proposition $2.2$, the only difference is the beginning of the proof, we apply the operator 
$\nabla$ to the equation \eqref{carlpb2} in the case \eqref{cond1th2} and 
the operators $\nabla\nabla\nabla$ in the case \eqref{cond2th2}. 
After that we exactly follow the steps $1$,$2$, and $3$ of the proof of Proposition $2.2$.

As a consequence, we also can derive the following observability inequality, whose proof is very classical (see also Proposition \ref{prop ine obs 2}):

\begin{prop}Under assumptions \eqref{cond1th2} and \eqref{reg1th2} or under assumptions \eqref{cond2th2} and \eqref{reg2th2}, 
 for every $\psi^0\in L^2(\Omega)^2$, the solution to System \eqref{carlpb2} satisfies
 \begin{equation}\label{ine obs 4}
 \begin{array}{c}
 \displaystyle\int_{\Omega} 
|\psi(0,x)|^2dx
\leqslant C_{obs}
\displaystyle\iint_{(0,T)\times\omega_1} e^{-2s_0\alpha}\xi^9
|\psi|^2dxdt,
\end{array} \end{equation}
where $s_0:=C(T^5+T^{10})$ and $C_{obs}:=e^{C(1+T+1/T^5)}$.
\end{prop}

To conclude, one can obtain the following controllability result:
\begin{prop}\label{prop cont reg2}
Consider the following system:
 \begin{equation}\label{primal 2x2 contr2}
 \left\{\begin{array}{ll}
\partial_tz_1=\Div(d_1\nabla z_1)+g_{11}\cdot\nabla   z_1
+g_{12}\cdot\nabla  z_2+a_{11}z_1+a_{12}z_2+\theta v_1&\mathrm{in}~ Q_T,\\

\partial_tz_2=\Div(d_2\nabla z_2)+g_{21}\cdot\nabla  z_1
+g_{22}\cdot\nabla  z_2+a_{21}z_1+a_{22}z_2+\theta v_2&\mathrm{in}~ Q_T,\\
       z_1=z_2=0&\mathrm{on}~\Sigma_T,\\
       z_1(0,\cdot)=y_1^0,~z_2(0,\cdot)=y_2^0&\mathrm{in}~\Omega.
        \end{array}
\right.
\end{equation} 
Under Conditions \eqref{cond1th2} and \eqref{reg1th2} or if $N=1$ and 
under Condition \eqref{cond2th2} and \eqref{reg2th2}, 
System \eqref{primal 2x2 contr2} is null controllable at time $T$, 
that is for every $y^0\in L^2(\Omega)^2$ 
there exists a control $v:=(v_1,v_2)\in L^2(Q_T)^2$ such that the solution $z$ 
to System \eqref{primal 2x2 contr2} satisfies $z(T)\equiv 0$ in $\Omega$. 
Moreover for every $K\in (0,1)$, we have   $e^{Ks_0\alpha^*}v\in \mathcal{X}^2$ where:
\begin{enumerate}
 \item[(i)]Under Conditions \eqref{cond1th2} and \eqref{reg2th2}, $\mathcal{X}$ is defined by
 \begin{equation}\label{defX1}
 \mathcal{X}:=L^2(0,T;H^4(\Omega)\cap H^1_0(\Omega))\cap H^2(0,T; L^2(\Omega)).
\end{equation}
 \item[(ii)]If $N=1$ and under Conditions \eqref{cond2th2} and \eqref{reg2th2}, $\mathcal{X}$ is defined by
 \begin{equation}\label{defX2}
 \mathcal{X}:=L^2(0,T;H^2(\Omega)\cap H^1_0(\Omega))\cap H^1(0,T; L^2(\Omega)).
\end{equation}
\end{enumerate}
Moreover in the both cases, we have the estimate
\begin{equation}\label{exp control2}
\|e^{Ks_0\alpha^*}v\|_{\mathcal{X}^2}\leqslant e^{C(1+T+1/T^5)}\|y^0\|_{L^2(\Omega)^2}.
\end{equation}
\end{prop}
One more time, the proof is very similar to the proof of Proposition $2.3$, notably and one can recover easily estimates  on the derivatives of order $1$ and $2$ in time of the control by using equation \eqref{carlpb2} verified by $\varphi$ and estimates similar to \eqref{estim H2 u 3}.

\subsection{Proof of Theorem \ref{th 2}}

  The proof is totally similar to the one of Theorem \ref{th 1}.
Let us assume that one of the two Conditions \eqref{cond1th2}  or \eqref{cond2th2} holds, and let us prove that System \eqref{syst21} is null controllable at time $T$ (which will imply the approximate controllability at time $T$).
Let $y^0\in L^2(\Omega)^2$. 
Using Proposition \ref{prop cont reg2},
System \eqref{primal 2x2 contr2}
 is null controllable at time $T$,  
more precisely there exists a control $v\in L^2(Q_T)^2$ such that the solution $z$ in $W(0,T)^2$ 
to System \eqref{primal 2x2 contr2} satisfies 
\begin{equation*}
z(T)\equiv 0 \mathrm{~~in~~} \Omega.  
\end{equation*}
Moreover  
\begin{equation}\label{exp control th2}
 e^{Ks_0\alpha^*}v\in \mathcal{X}^2.
\end{equation}
Let us remind that $\theta$ was defined in \eqref{theta2}.
Let $(\widehat{z}_1,\widehat{z}_2,\widehat{v})$ be defined by
\begin{equation*}
\left(\begin{array}{c}\widehat{z}_1\\\widehat{z}_2\\\widehat{v}\end{array}\right)
 :=\mathcal{M}\left(\begin{array}{c}\theta v_1\\\theta v_2\end{array}\right),
\end{equation*}
with $\mathcal M:\mathcal{X}^2\rightarrow L^2(Q_T)^2\times L^2(Q_T)$ given as the adjoint of $\mathcal M^*$ defined in \eqref{M*21} and \eqref{M*22}. 
Thanks to the definition of $\theta$ given in \eqref{theta2}, the fact  the coefficients of $\mathcal M$ are necessarily at least in $L^\infty((a,b)\times \omega_0)$, 
the definition of $\mathcal{X}$ given in \eqref{defX1}-\eqref{defX2} and the fact that $\mathcal M$ 
is of order $1$ in time and $2$ in space under Condition \eqref{cond1th2} and 
is of order $2$ in time, $4$ in space and $1-2$ in crossed time-space 
(which is an interpolation space between $L^2((0,T),H^4(\Omega))$ and $H^2((0,T),L^2(\Omega))$ thanks to \cite[13.2, P. 96]{lions1968problemes}) under Condition \eqref{cond2th2},
we obtain $(\widehat{z}_1,\widehat{z}_2,\widehat{v})\in L^2(Q_T)^2\times L^2(Q_T)$.
Moreover, using \eqref{exp control th2}, we remark that $(\widehat{z}_1,\widehat{z}_2,\widehat{v})$ 
is a solution   to the control problem
\begin{equation}\label{probleme ramene a tout l espace th2}
 \left\{\begin{array}{ll}
       \partial_t\widehat{z}=\Div(D\nabla\widehat{z})+G\cdot\nabla \widehat{z}+A\widehat{z}+B\widehat{v}+\theta v&\mathrm{in}~ Q_T,\\
       \widehat{z}=0&\mathrm{on}~\Sigma_T,\\
       \widehat{z}(0,\cdot)=\widehat{z}(T,\cdot)=0&\mathrm{in}~\Omega.
        \end{array}
\right.
\end{equation}
Finally, $\widehat{z}\in W(0,T)^2$ thanks to the usual parabolic regularity. 
Thus the pair $(y,u):=(z-\widehat{z},-\widehat{v})$ is a solution to System \eqref{syst21} 
in $W(0,T)^2\times L^2(Q_T)$ and satisfies 
\begin{equation*}
y(T)\equiv 0 \mathrm{~~in~~} \Omega,
\end{equation*}
which concludes the proof of Theorem \ref{th 2}.

\section*{Acknowledgements}
The authors would like to thank the anonymous referees for their numerous comments and corrections that greatly contributed to improving the final version of the paper.



\end{document}